\documentclass[final,5p,times,twocolumn,authoryear]{elsarticle}
\usepackage{amsmath}
\usepackage{amssymb}
\usepackage[round]{natbib}            
\usepackage{graphicx}          
\usepackage{mathdots}
\usepackage{textcomp}
\usepackage{comment}
\usepackage{stfloats}
\usepackage[shortlabels]{enumitem}

\begin{document}

\begin{frontmatter}

\title{Lyapunov modal analysis and participation factors with applications
\\ to small-signal stability of power systems}

\author[First]{Alexey B. Iskakov}\ead{iskalexey@gmail.ru}\corref{*}
\author[First]{Igor B. Yadykin}\ead{jad@ipu.ru}
\cortext[*]{Corresponding author: A.~Iskakov, Laboratory 19, V.A.~Trapeznikov Institute of Control Sciences RAS,
Ul. Profsoyuznaya 65, Moscow 117997, Russia. Tel.: +7 926 859 2632; +7 495 334 9030. Fax: +7 495 334 9340. E-mail:
 iskalexey@gmail.ru, isk\_alex@mail.ru}
\address[First]{V.A.~Trapeznikov Institute of Control Sciences of RAS, Profsoyuznaya 65, Moscow 117997, Russia \\
(e-mail: iskalexey@gmail.ru, jad@ipu.ru)}
          

\begin{abstract}                          
When random disturbances are regularly introduced into a dynamical system over time, its small-signal stability is determined by the energy of perturbations accumulated in the system. 
To analyze this perturbation energy, this paper proposes a novel physically motivated Lyapunov modal analysis (LMA) framework, which combines selective modal analysis with the spectral decompositions of specially chosen Lyapunov functions. This approach allows  the modal interactions in dynamical systems to be characterized and estimated in connection with specific state variables. Conventional participation factors characterize the relative contribution of the system modes and state variables to the evolution of states and modes, respectively. In contrast, the proposed Lyapunov participation factors characterize similar contributions to corresponding Lyapunov functions, which determine the integral energy associated with the states and modes on an infinite or finite time interval. This allows the estimation of modal interactions in terms of total energy produced by their mutual actions over time. Using a two-area four-machines power system, we demonstrate that LMA reliably identifies resonant modal interactions, merging of modes, and loss of stability, even for a linear model, and associates them with certain state variables. The Lyapunov participation factors corresponding to the selected part of the system spectrum can be calculated independently and serve as a basis for rapid  real-time calculations of critical mode behaviors in large-scale dynamical systems. 
\end{abstract}

\end{frontmatter}

\section{Introduction}

Developing modern 
smart grid and microgrid technologies are a priority in the advancement of electric power systems (EPSs). A critical requirement for the introduction of these technologies is increasing the reliability of EPS and the ability to monitor and control its stability in real time \citep{hager:2014}. 
In a large EPS, weakly stable modes usually occur in groups, leading to resonance interaction problems and the appearance of dangerous 
low-frequency oscillations. Such oscillations may occur within the power facility, regional power grid, or global power system \citep{pal:2012}. Quite frequently, these oscillations establish critical limitations for maximum transfer capability in the main power transmission lines and
lead to the occurrence of voltage collapse and cascading failures \citep{weber:2016}. The loss of stability is accompanied by the accumulation of energy in the low-frequency oscillations that causes a resonant reaction in the system. Therefore, the small-signal stability analysis of modern EPSs requires an accurate estimation and prediction of the resonant interactions of weakly stable system modes with reference to specific state variables.

In the conventional modal analysis of linear systems, modal interactions are not taken into account directly. In nonlinear versions of modal analysis, modal interactions are considered in the second and higher order terms of the Taylor expansion of the system dynamics. In general, this approach imposes demanding requirements on model accuracy and has a high computation cost. More importantly, however,
conventional modal analysis 
estimates small-signal stability with respect to some initial disturbance in the system  
and modal interactions are considered in terms of instantaneous dynamics.
However, when random disturbances are introduced into the system regularly, the time factor becomes critical for the small-signal stability.
In this case, the system stability is determined not so much by the instantaneous dynamics of a single perturbation, but rather by the energy of perturbations accumulated in the system over time. This perturbation energy can be estimated using the spectral decompositions of Lyapunov functions proposed in \citep{yad:10,yad:17}. 
For this purpose, this paper proposes a Lyapunov modal analysis (LMA) framework 
that combines selective modal analysis with the spectral decompositions of specially selected Lyapunov functions.
This approach allows the estimation of modal interactions based on the total energy produced in a system by mutual action of modes over time rather than in terms of their instantaneous dynamics. 
It is conceptually different from considering second-order and higher terms in the instantaneous dynamics of the nonlinear model.

\subsection{Literature review}

Modal analysis is one of the most popular methods for studying the small-signal stability of dynamical systems. 
A {\it selective modal analysis (SMA)} framework proposed in \citep{perez:82,verghese:82} allowed an accurate identification of the elements of the system structure associated with specific eigenmodes based on the so-called {\it participation factors (PFs)}. For linear time-invariant systems, PFs have been defined as the relative contributions of state variables to the evolution of system modes, or as the relative contributions of system modes to the evolution of states. Subsequently, the concept of PFs gained widespread use in power engineering and other applications for analyzing stability \citep{verghese:82, song:2019}, reducing dynamic models \citep{chow:2013}, determining the optimal placement of sensors and stabilizers \citep{singh:2010}, and solving clustering problems \citep{genk:2005}. The interpretation of PFs has been expanded in terms of the sensitivity of eigenvalues \citep{pagola:89}, modal controllability and observability \citep{hamdan:89}, and modal mobility \citep{hamdan:2010}. This original interpretation of PFs involved specially selected initial conditions, and it was observed that this assumption could lead to counterintuitive results \citep{hash:09}; therefore, an alternative method of averaging over an uncertain set of initial system conditions was proposed. Accordingly, the original definition of PFs was retained for {\it the “mode-in-state” PFs}; however, a new definition was proposed for the {\it “state-in-mode” PFs (or SIMPFs)}. Subsequently, similar SIMPF concepts were considered for dynamical nonlinear systems \citep{hamzi:2014} and for systems described by algebraic equations such as power flow equations \citep{song:2019}. 

Attempts to account for {\it nonlinear effects} and {\it intermodal interactions} within the framework of modal analysis developed mainly in two directions. The model-based approach is associated with taking into account the higher-order terms of the Taylor expansion in the system approximation and using normal Poincar\'e forms \citep{vittal:91, hamzi:2014, tian:2018}. A study in \citep{sanchez:2005} showed that accounting for these terms becomes significant, for example, when studying inter-area oscillations in stressed power systems following large disturbances. These methods, however, generally require solving a highly nonlinear numerical problem using computationally expensive algorithms. Another approach involves estimating the PFs directly from measurements. This approach, for example, can be based on extended  dynamic mode decomposition \citep{williams:2015} and Koopman mode decomposition \citep{netto:2019}. The performance of these methods still requires careful verification in practical applications. 

Another conceptual method in stability analysis is associated with the use of Lyapunov equations \citep{dahleh:2011}. {\it Lyapunov stability analysis} is based on a positive definite function of a system state $V({\boldsymbol{x}}) = \boldsymbol{x}^T {\bf P} \, \boldsymbol{x}$, where the positive definite matrix ${\bf P} = {\bf P}^T > 0$ is a solution for the corresponding Lyapunov equation which is called {\it the Gramian}. For linear time-invariant systems, Lyapunov functions can be associated with the integrated energy of the input or output signal. The Gramians of controllability and observability are commonly used in connection to this. In general, the {\it observability Gramian} characterizes the system stability in terms of its output energy limit while the {\it controllability Gramian} does so in terms of its asymptotic sensitivity to the random input disturbances. For a stable linear system, the Gramians are closely related to the squared $H_2$ norm of its transfer or impulse response functions. The physical interpretation of these values is that they characterize energy amplifications in the system averaged over time or frequency. The energy-based interpretation of Gramians generally holds true for time-varying linear systems if the exponential expressions $e^{{\bf A}t}$ are replaced with a fundamental solution ${\bf\Phi}(0, t)$ of the homogeneous equation $\dot{\boldsymbol{x}} = {\bf A}(t) \, \boldsymbol{x}$ \citep{shokoohi:83, verriest:83}. The concept of Gramians was further generalized and interpreted for deterministic bilinear and stochastic linear systems as {\it energy functionals} \citep{gray-mesko:98, benner-damm:11}.

Lyapunov stability analysis was applied to assess the stability of electric power systems in \citep{pavella:2012,chiang:2011}. 
The spectral properties of Gramians and energy functionals were effectively used for {\it model order reduction} techniques. These techniques include methods for balanced truncation \citep{moore:81}, the use of cross-Gramians \citep{fernando:84} and their various modifications (see the review in \citet{baur:2014}). \citet{antoulas:2005} 
obtained singular expansions for infinite Gramians of controllability and observability based on diagonalization of the dynamics matrix. A more general form of the {\it spectral decompositions of Lyapunov functions} into  components corresponding to the individual eigenvalues of the system or their pairwise combinations was proposed in \citep{yad:10, yad:14, yad:17, zubov:17}. Each eigen-part was denoted a {\it sub-Gramian}. These allow estimations of the interactions between eigenmodes and were applied to the stability analysis of power system in \citep{yad-cep:16}.

\subsection{Main contribution}

The objective of this paper is to offer a novel framework for {\it Lyapunov modal analysis} that combines two approaches, selective modal analysis and the spectral decompositions of Lyapunov functions, for stability assessment. For this purpose, we propose the concept of {\it Lyapunov participation factors}, which characterize the relative contribution of system modes $\boldsymbol{z}$ and state variables $\boldsymbol{x}$, not for the evolution of states and modes, respectively, but for the corresponding Lyapunov functions, i.e., to the quantities $\boldsymbol{x}^T{\bf P}_x\boldsymbol{x}$ or $\boldsymbol{z}^*{\bf P}_z \boldsymbol{z}$. The matrices ${\bf P}_z$ and ${\bf P}_x$ here are the solutions of the Lyapunov equations, which are chosen such that their solutions measure the integrated energy associated with a particular eigen-mode or state variable. The corresponding Lyapunov functions are denoted {\it Lyapunov energies}. These values are important when analyzing the stability of the system because they do not reflect the instantaneous values of the states or signals, but the variation in their energies over a time interval (i.e., their energy  gains in the system). In terms of mechanics, Lyapunov energy corresponds more to {\it Hamilton's action}, i.e., the energy integrated over time \citep{landau:78}, than to the energy itself.

The key question is how to define the energy of states and modes. For the former, we follow the definition of {\it stored energy} proposed by \citet{MacFarlane:69} for electromechanical systems, which is a quadratic function $\frac{1}{2}\,\boldsymbol{x}^T\boldsymbol{x}$ of state variables $\boldsymbol{x}$ after suitable scaling. 
After this, the Lyapunov energy of the $k$-th state variable $x_k$ can be defined as
\begin{equation} \nonumber
E_{x_k} =\int\limits_0^{t} x^2_k (\tau) d\tau  \ . 
\end{equation}
The definition of modal energy is less obvious. The modal energy of the $i$th mode was defined in \citep{hamdan:86} as $\boldsymbol{x}^T \boldsymbol{u}_i (\boldsymbol{\varv}_i)^T \boldsymbol{x}$, where $\boldsymbol{u}_i $ and $\boldsymbol{\varv}_i$ are the normalized right and left (column) eigenvectors of the $i$th mode. Unfortunately, this definition lacks a physical meaning; using it results in
modal energies that are negative and unlimited in magnitude even in a very simple system (see \citet{isk:19}).

Therefore, two alternative definitions of modal energy are investigated in this study. First, $i$-th mode energy can be defined as 
\[
 (\boldsymbol{x})^T (\boldsymbol{\varv}^T_i)^* (\boldsymbol{u}_i)^*  \boldsymbol{u}_i \, \boldsymbol{\varv}^T_i \boldsymbol{x} = |\boldsymbol{u}_i|^2 |z_i|^2 \ ,
\]
where $z_i = \boldsymbol{\varv}^T_i \boldsymbol{x}$ is the $i$th component of the system mode vector $\boldsymbol{z}$.
Then, choosing the normalization of eigenvectors so that $|\boldsymbol{u}_i|^2=1$, the Lyapunov energy of the $i$-th mode can be represented as 
\[
E_{z_i} =\int\limits^{t}_0 |z_i(\tau)|^2  d\tau  \ , 
\]
by analogy with $E_{x_k}$.
We show that this definition leads in practice to the “state-in-mode” PFs defined in \citep{hash:09} for real eigenvalues and corrected in \citep{konoval:17, isk:19} for the case of complex eigenvalues.
Alternatively, the energy of the $i$th mode can be defined as a modal contribution of the $i$th mode to the Lyapunov energy of states (i.e., into $\sum_k E_{x_k}$) based on the spectral decompositions proposed in \citep{yad:10, yad:17}. With this definition, modal contributions of some modes can be negative if there are other modes with a sufficiently large amplitude and negative correlation with a given mode. This correlation between modes is determined by both the dynamic properties of the system and its current state. We introduce the idea of {\it Lyapunov modal interaction energies and factors} that characterize pairwise modal interactions in terms of the Lyapunov energy they produce in different state variables, based on this second definition of modal energy. We also offer two indicators for selecting the state variables that are the most (i) sensitive for identifying a specific modal interaction and (ii) influential for dampening this interaction.

The definitions of Lyapunov energies for states and modes allow us to introduce the corresponding concepts of {\it Lyapunov participation factors}. We examine the theoretical properties of these indicators and use a simulation of a power system \citep{kundur:94} to show that they can reliably identify resonant modal interactions, mode merging, and the loss of stability, and associate these events with certain state variables.

Unlike works on nonlinear modal analysis \citep{vittal:91,hamzi:2014,tian:2018}, this paper proposes a new principle for evaluating modal indicators and modal interactions. This principle is not based on the instantaneous dynamics of variables, but on variation of their energy over a time interval. It allows the identification of low-frequency oscillations dangerous for small-signal stability and the detection of the effect of energy accumulation in these oscillations. Unlike works on spectral decompositions of Lyapunov functions \citep{yad:14,yad-cep:16,yad:17,zubov:17}, the proposed method allows the association of these decompositions with specific state variables and their application to problems of modal analysis. The idea of combining selective modal analysis and Lyapunov spectral expansions has been already mentioned in \citep{vassil:17}. This paper presents a general framework for its implementation.

\subsection{Organization of the paper}

Preliminary information on PFs and sub-Gramians is briefly summarized in Section 2. 
Section 3 introduces the Lyapunov energies of states and modes and the corresponding Lyapunov PFs.
Lyapunov modal interaction analysis and the corresponding pair PFs are first mentioned and defined in Section 4. 
In Section 5, some characteristic properties of Lyapunov PFs are established. 
The numerical experiment that demonstrates the potential advantages of Lyapunov modal analysis is provided in Section 6.
Section 7 contains the conclusions drawn from this work.

\section{Theoretical background} 

\subsection{Participation Factors}

In this subsection, we recall the definition and some properties of the participation factors given in 
\citep{perez:82, pagola:89, garofalo:02}.
Consider an autonomous linear time-invariant system 
\begin{equation} \label{system0}
\dot{\boldsymbol{x}} = {\bf A} {\boldsymbol{x}}(t) \ ,
\end{equation}
where $\boldsymbol{x} \in \mathbb{R}^n$ is a system state vector and ${\bf A}\in \mathbb{R}^{n\times n}$ is a real matrix with a simple spectrum 
that can be represented as
\begin{equation} \label{A=ULV}
{\bf A} = {\bf U}{\bf \Lambda} {\bf V} = ({\boldsymbol{u}}_1 {\boldsymbol{u}}_2 \dots {\boldsymbol{u}}_n) 
\begin{pmatrix}
    \lambda_1       & 0 &  \dots & 0 \\
    0      & \lambda_2 &  \dots & 0 \\
    \vdots & \vdots  & \ddots  & \vdots \\
    0      & 0  & \dots & \lambda_n
\end{pmatrix}
\begin{pmatrix}
  {\boldsymbol{\varv}}^T_1 \\
  {\boldsymbol{\varv}}^T_2 \\
   \vdots  \\
   {\boldsymbol{\varv}}^T_n 
\end{pmatrix},
\end{equation}
where ${\bf U} {\bf V} = {\bf V} {\bf U} = {\bf I}$ and $(\cdot)^T$ is the transpose operation. 
Let $u_i^k$ and $\varv_i^l$ be the $k$-th and $l$-th components of the eigenvectors 
${\boldsymbol{u}}_i$ and ${\boldsymbol{\varv}}_i$, respectively. Then
\begin{equation} \label{PF-def}
p_{ki} = u_i^k \varv_i^k \ \ \text{and} \ \  p_{kil} = u_i^k \varv_i^l
\end{equation}
are called {\it participation factors} (PFs) and {\it generalized participations}, respectively \citep{pagola:89}.
The PF $p_{ki}$ $``$weights$"$ the participation of the $i$-th mode in the $k$-th state variable and 
was interpreted in \citep{hash:09} as the {\it mode-in-state PF}. 
We then recall the following two important properties of generalized participations from \citep{garofalo:02}. \\

\noindent {\bf Property $1$}. The generalized participation $p_{kil}$ is the 
sensitivity of the $i$-th eigenvalue $\lambda_i$ with respect to the element $a_{lk}$ of ${\bf A}$, i.e.,
\begin{equation}\label{PF-sensetivities}
p_{kil}=\frac{\partial \lambda_i}{\partial a_{lk}} \ , \ \ p_{ki} = \frac{\partial \lambda_i}{\partial a_{kk}}  \ .
\end{equation}

The residue matrices ${\bf R}_i$ are defined as the coefficients in the expansion of the resolvent matrix ${\bf A}$ \citep{garofalo:02}:
\begin{equation} \label{residues}
(s {\bf I} - {\bf A})^{-1} = \frac{{\bf R}_1}{s-\lambda_1} + \frac{{\bf R}_2}{s-\lambda_2} + \dots + \frac{{\bf R}_n}{s-\lambda_n} \ .
\end{equation}

\noindent {\bf Property $2$}. The generalized participations $p_{kil}$ 
are the coefficients of the corresponding residue matrix ${\bf R}_i$, i.e.,
\begin{equation} \label{prop-2}
p_{kil} = {\bf e}_k^T \ {\bf R}_i \  {\bf e}_l \ , \ \  {\bf R}_i = \sum_{k,l} p_{kil} \ {\bf e}_k \cdot {\bf e}_l^T \ , 
\end{equation}
where ${\bf e}_k$ and ${\bf e}_l$ are the $k$-th and $l$-th columns of the identity matrix.

By applying 
to \eqref{system0} the diagonalizing transformation 
\begin{equation}  \label{z-vector}
\boldsymbol{z}= {\bf V} \boldsymbol{x} \ , 
\end{equation}
the evolution of the system mode vector $\boldsymbol{z}(t)$ can be specified as 
\begin{equation} \label{z-evolution}
\dot{\boldsymbol{z}}(t) = {\bf V} \,{\bf A}\boldsymbol{x}(t) = {\bf \Lambda} {\bf V} \boldsymbol{x}(t) = {\bf \Lambda} \boldsymbol{z}(t) \ .
\end{equation}
For the analysis of $\boldsymbol{z}(t)$ 
\citet{hash:09} proposed the use of {\it state-in-mode PFs} 
which, taking into account the correction in \citep{konoval:17}, were defined as
\begin{equation} \label{simpf}
{\pi }_{ki}=\frac{(\varv^k_i)^* \, \varv^k_i}{{{\left({\boldsymbol{\varv}}_i\right)}^* \, \boldsymbol{\varv}}_i }\ \ .\ 
\end{equation}

\subsection{Gramians and sub-Gramians}

Here, we recall the definition and some properties of the Gramians and sub-Gramians 
from \citep{yad:10, yad:14, yad:17}. 
A Gramian is a positive definite solution ${\bf P}= {\bf P}^*>0$ of the following matrix Lyapunov equation:
\begin{equation} \label{lyapunov}
  {\bf A}^* \, {\bf P} + {\bf P} \, {\bf A} = - \ {\bf Q} \ , \ \ {\bf Q} = {\bf Q}^* > 0 \ ,
\end{equation}
where $(\cdot)^*$ denotes the conjugate transpose of a matrix. 
For simplicity, we also assume that the matrix ${\bf A}$ has a simple spectrum. 
In this case, the solution of (\ref{lyapunov}) can be written as follows \citep{yad:17}: 
\begin{equation} \label{pair-subgram}
{\bf P} = - \sum\limits_{i,j=1}^n \frac{{\bf R}^*_i \, {\bf Q} \, {\bf R}_j}{\lambda_i^*+\lambda_j} \ ,
\end{equation}
where ${\bf R}_i$ and ${\bf R}_j$ are the residue matrices defined by (\ref{residues}) 
corresponding to $\lambda = \lambda_i$ and $\lambda = \lambda_j$, respectively. 

It follows from (\ref{residues}) that
\begin{equation} \nonumber
\sum\limits_{j=1}^n \frac{{\bf R}_j}{-\lambda^*_i - \lambda_j} = (-\lambda^*_i {\bf I} - {\bf A})^{-1} \ .
\end{equation}
Substituting this in (\ref{pair-subgram}), we obtain another form of the spectral decomposition as follows: 
\begin{equation} \label{subgram}
{ \bf P} = - \sum^n_{i=1} {\bf R}^*_i \ {\bf Q} \ (\lambda^*_i {\bf I}+{\bf A})^{-1} 
= - \sum^n_{j=1} (\lambda_i {\bf I} + {\bf A}^*)^{-1} {\bf Q} \ {\bf R}_j \ , 
\end{equation}
The following Hermitian parts of matrices in the spectral decompositions (\ref{pair-subgram})
and (\ref{subgram}) 
 \begin{equation} \label{subgram+}
 \widetilde{\bf P}_i = - \left\{{\bf R}^*_i \ {\bf Q} \ (\lambda^*_i {\bf I} + {\bf A})^{-1} \right\}_H  , \ \ 
 {\bf P}_{ij}=-\left\{\cfrac{{\bf R}^*_i \ {\bf Q} \ {\bf R}_j}{\lambda^*_i + \lambda_j }\right\}_H \ .
\end{equation}
have been denoted {\it sub-Gramians}. Here, $\{\cdot\}_H$ represents the Hermitian part of a matrix, 
i.e., $\left\{{\bf M} \right\}_H = \frac{1}{2}\left({\bf M} + {\bf M}^*\right)$. 
Using sub-Gramians (SGs), the decompositions (\ref{pair-subgram}) and (\ref{subgram})  can be written as 
\begin{equation} \label{decomp}
{\bf P} = \sum_{i=1}^n \widetilde{\bf P}_i = \sum_{i,j=1}^n {\bf P}_{ij} \ , \ \   \widetilde{\bf P}_i = \sum_{j=1}^n {\bf P}_{ij} \ .
\end{equation}
The SGs $\widetilde{\bf P}_i$ and ${\bf P}_{ij}$ characterize the contribution of eigenmodes or their pairs 
to the asymptotic variation of perturbation energy in the system defined by the Gramian ${\bf P}$.

The SGs \eqref{subgram+} also satisfy the Lyapunov equations, as follows:
\begin{gather}
  {\bf A}^* \, \widetilde{\bf P}_i + \widetilde{\bf P}_i \, {\bf A} = - \ \frac{1}{2}\Big({\bf R}^*_i {\bf Q}+{\bf Q} \,{\bf R}_i\Big)  \ , \nonumber \\
  {\bf A}^* \, {\bf P}_{ij} + {\bf P}_{ij} \, {\bf A} = - \ \frac{1}{2}\Big({\bf R}^*_i {\bf Q}\, {\bf R}_j + {\bf R}^*_j {\bf Q} \, {\bf R}_i\Big)  \ . \label{lyap-subgram}
\end{gather}
These equations can be verified by substituting \eqref{subgram+} directly into them.\\

\subsection{Relationship between participation factors (PFs) and sub-Gramians (SGs)}

Given all of the above, the relationship between PFs and SGs can now be established. Although 
they are clearly different in terms of their conceptual definitions,
it is easy to see that both quantities are calculated using matrix residues. Substituting expression (\ref{prop-2}) for the matrix
residues ${\bf R}_i$ and ${\bf R}_j$ in (\ref{subgram+}), we obtain the following:
\begin{equation} \label{subgram-PF}
 \widetilde{\bf P}_i = -\sum\limits^n_{k,l=1} \left\{p^*_{lik} \ {\bf e}_k  {\bf e}^T_l \,
 {\bf Q} \, (\lambda^*_i {\bf I} +{\bf A})^{-1}\right\}_H \ ,
\end{equation}
\begin{equation}  \label{pair-subgram-PF}
{\bf P}_{ij} = - \sum\limits^n_{k,l,m,r=1} \left\{\frac{p^*_{lik} \ p_{rjm} \ {\bf e}_k  {\bf e}^T_l  {\bf Q} \ {\bf e}_r {\bf e}^T_m} {\lambda^*_i + \lambda_j} \right\}_H .
\end{equation}
These new representations allow the calculation of SGs via PFs. The SG $ \widetilde{\bf P}_i$
is the sum of terms proportional to the PFs corresponding to
$i$\nobreakdash-th eigenmode. The pairwise SG ${\bf P}_{ij}$ is a quadratic form of
PFs corresponding to the $i$\nobreakdash-th and $j$\nobreakdash-th eigenmodes. Formulas
similar to (\ref{subgram-PF}) were obtained in (\cite{antoulas:2005}, p.149).

Substituting (\ref{PF-def}) into (\ref{subgram-PF}) and (\ref{pair-subgram-PF}), we obtain the following:

\begin{equation} \label{subgram-eigen}
 \widetilde{\bf P}_i = - \left\{ ({\boldsymbol{\varv}}^T_i)^* {\boldsymbol{u}}^*_i 
 \, {\bf Q} \ (\lambda^*_i {\bf I} +{\bf A})^{-1}\right\}_H \ ,
\end{equation}
\begin{equation} \label{pair-subgram-eigen}
{\bf P}_{ij} = - \left\{ \frac{ ({\boldsymbol{\varv}}^T_i)^* {\boldsymbol{u}}^*_i 
\, {\bf Q} \ {\boldsymbol{u}}_j \, {\boldsymbol{\varv}}^T_j} 
{\lambda^*_i + \lambda_j}\right\}_H \ .
\end{equation}
These equations allow the calculation of SGs using eigenvectors that correspond to them when ${\bf A}$ has a simple spectrum. They demonstrate that the calculation of the individual SGs does not require knowledge of the entire spectrum, but only of the eigenvalues and eigenvectors corresponding to the particular SG.

Substituting (\ref{PF-sensetivities}) into (\ref{subgram-PF}) and (\ref{pair-subgram-PF}), we obtain
\begin{equation} \nonumber
\begin{array}{c}
 \widetilde{\bf P}_i = - \sum\limits^n_{k,l=1} \left\{  \,\left(\frac{\partial\lambda_i}{\partial a_{kl}}\right)^*
 {\bf e}_k  {\bf e}^T_l \, {\bf Q} \ (\lambda^*_i {\bf I} +{\bf A})^{-1} \right\}_H \ , \\
{\bf P}_{ij} = -  \sum\limits^n_{k,l,m,r=1} \left\{ \frac{1}{\lambda^*_i + \lambda_j} \left(\frac{\partial\lambda_i}{\partial a_{kl}}\right)^*
\frac{\partial\lambda_j}{\partial a_{mr}} \, {\bf e}_k {\bf e}^T_l \, {\bf Q} \ {\bf e}_r {\bf e}^T_m \right\}_H .
\end{array}
\end{equation}
These equations allow the calculation of SGs using the sensitivities of the corresponding eigenvalues with respect to the elements of matrix 
${\bf A}$.

\section{Lyapunov energies and participation factors}

In this section, we introduce new indicators for selective modal analysis.
Unlike conventional PFs, the proposed Lyapunov PFs characterize the relative contribution 
of the system modes~$\boldsymbol{z}$ and state variables~$\boldsymbol{x} $ 
not to the evolution of states and modes, respectively, 
but to the corresponding {\it Lyapunov energies}, i.e., to the quantities $\boldsymbol{x}^T{\bf P_{\boldsymbol x}}\, \boldsymbol{x}$ or
$\boldsymbol{z}^*{\bf P_{\boldsymbol z}}\, \boldsymbol{z}$, where the {\it Gramians} 
${\bf P_{\boldsymbol x}}$ and ${\bf P_{\boldsymbol z}}$ are  
the solutions of the corresponding Lyapunov equations.
For this purpose, we first consider the concept of {\it Lyapunov energies}, which characterize the squared
state variables and system modes integrated over time. 

\subsection{Lyapunov energies of states and modes} 

Consider a linear dynamical system of the form:
\begin{equation} \label{system}
\dot{\boldsymbol{x}} = {\bf A} {\boldsymbol{x}} \ , \ \ {\boldsymbol{y}} = {\bf C}{\boldsymbol{x}} \ , 
\end{equation}
where $\boldsymbol{x}, \boldsymbol{y}\in \mathbb{R}^n$ are the state and output vectors,
$\boldsymbol{\mathrm{C}}\in \mathbb{R}^{n\times n}$ is an output matrix, and $\boldsymbol{\mathrm{A}}\in \mathbb{R}^{n\times n}$ is a stable state matrix, 
which has $n$ distinct eigenvalues $({\lambda }_1,{\lambda }_2,\ldots ,{\lambda }_n)$ and 
can be represented through its eigenvectors by~\eqref{A=ULV}.
In further analysis, we also assume that the eigenvectors associated with real eigenvalues are real, and the eigenvectors associated with a complex conjugate pair of nonreal eigenvalues ${\lambda }_i,{\lambda }^*_i$ are complex conjugates: 
\begin{equation} \label{ev-assumption}
\begin{array}{c}
\rm{Re}\{\boldsymbol{u}_{\it{i*}}\} = \rm{Re}\{\boldsymbol{u}_{\it{i}}\}, \ \ \rm{Im}\{\boldsymbol{u}_{\it{i*}}\} = -\rm{Im}\{\boldsymbol{u}_{\it{i}}\} , \\
\rm{Re}\{\boldsymbol{\varv}_{\it{i*}}\} = \rm{Re}\{\boldsymbol{\varv}_{\it{i}}\}, \ \ \rm{Im}\{\boldsymbol{\varv}_{\it{i*}}\} = -\rm{Im}\{\boldsymbol{\varv}_{\it{i}}\} ,
\end{array}
\end{equation}
where $i*$ is the index of the eigenvector associated with ${\lambda }^*_i$. 
Taking into account the normalization conditions \eqref{ev-assumption}, the eigenvectors in the representation \eqref{A=ULV} are not uniquely determined, but up to an arbitrary constant using which for any $i$, one can multiply the left eigenvector $\boldsymbol{\varv}_i$ and divide the right eigenvector $\boldsymbol{u}_i$. 
Therefore, to avoid this ambiguity, we impose an additional normalization condition
\begin{equation} \label{ev-assumption-2}
\forall i\, : \ \ \ |\boldsymbol{u}_i|^2 = 1 \ .
\end{equation}

The integrated output signal energy in the stable system \eqref{system} generated by
the initial state ${\boldsymbol{x}}_0$ can be defined as
\begin{equation} \label{energy-integral}
E = \int^{\infty}_0 {\boldsymbol{x}}^T_0 \, e^{{\bf A}^*t} \ {\bf C}^T {\bf C} \ e^{{\bf A}t} {\boldsymbol{x}}_0 \ dt 
= {\boldsymbol{x}}^T_0 \, {\bf P} \, {\boldsymbol{x}}_0 \ ,
\end{equation}
where the Gramian ${\bf P} = \int_0^{\infty} e^{{\bf A}^*t} \, {\bf C}^T {\bf C} \ e^{{\bf A}t} dt $ is the solution
of the Lyapunov equation ${\bf A}^*{\bf P} + {\bf P} \ {\bf A} = - \ {\bf C}^T {\bf C}$. 

In the matrix ${\bf C}$, let us take a unit vector with a unit $k$-state
component and other components equal to zero: 
\begin{equation} \nonumber
{\bf C} = {\bf e}^T_k =(0,\dots,0,1_k,0,\dots,0) \ .
\end{equation}
Then, $y = x_k$ is a $k$-th component of the state vector ${\boldsymbol{x}}$, 
and the expression \eqref{energy-integral} characterizes the integrated energy produced
in the $k$-th state. 
Therefore, we define the {\it Lyapunov energy produced in the $k$-th state} as 
\begin{gather}
E_{x_k} = {\boldsymbol{x}}^T {\bf P}_{x_k} {\boldsymbol{x}} = \int^{\infty}_0 x^2_k(t)\,dt \ , \ \ \text{where} \nonumber \\
 {\bf A}^* \, {\bf P}_{x_k} + {\bf P}_{x_k} \, {\bf A} = - \ {\bf e}_k {\bf e}^T_k \ . \label{Lener-k-state} 
\end{gather}
The total Lyapunov energy of all the state variables is 
\begin{equation} \label{LE-all-states}
E_x = \sum^n_{k=1} E_{x_k} = {\boldsymbol{x}}^T {\bf P}_{\boldsymbol{x}} \, {\boldsymbol{x}} \ , \ \text{where} \ \ {\bf P}_{\boldsymbol{x}} = \sum^n_{k=1} {\bf P}_{x_k} \ ,
\end{equation}
which is the solution of the Lyapunov equation 
\begin{equation} \label{Lyap-states}
{\bf A}^* \, {\bf P}_{\boldsymbol{x}} + {\bf P}_{\boldsymbol{x}} \, {\bf A} = - \ {\bf I} \ . 
\end{equation}

By applying 
the diagonalizing transformation \eqref{z-vector} to the state vector~$\boldsymbol{x}$,
we obtain the system mode vector~$\boldsymbol{z} = {\bf V} \boldsymbol{x} \in \mathbb{C}^n$, 
for which the equation of the dynamic \eqref{system} takes the diagonal form
$\dot{\boldsymbol{z}} = {\bf \Lambda} \boldsymbol{z}$.
By analogy with \eqref{Lener-k-state}, we define the {\it Lyapunov energy of the $i$-th mode} as 
\begin{gather}
E_{z_i} = \int^{\infty}_0 z^*_i(t) \, z_i(t)\,dt = \boldsymbol{z}^* \, {\bf P}_{z_i} \, \boldsymbol{z} \ , \ \ \text{where} \nonumber \\
{\bf \Lambda}^* \, {\bf P}_{z_i} + {\bf P}_{z_i} \, {\bf \Lambda} = - \ {\bf e}_i {\bf e}^T_i \ , \label{Lener-i-mode+}
\end{gather}
where $z_i$ is the $i$-th component of the mode vector~$\boldsymbol{z} \in \mathbb{C}^n$. 
The total Lyapunov energy of the mode vector~$\boldsymbol{z}$ is 
\begin{equation} \label{LE-all-modes+}
E_z = \int^{\infty}_0 \boldsymbol{z}^*(t) \, \boldsymbol{z}(t)\,dt =
\sum^n_{i=1} E_{z_i} = {\boldsymbol{z}}^* {\bf P}_{\boldsymbol{z}} \, {\boldsymbol{z}} \ , 
\end{equation}
where ${\bf P}_{\boldsymbol{z}} = \sum^n_{i=1} {\bf P}_{z_i}$
is the solution of the Lyapunov equation 
\begin{equation} \label{Lyap-modes}
{\bf \Lambda}^* \, {\bf P}_{\boldsymbol{z}} + {\bf P}_{\boldsymbol{z}} \, {\bf \Lambda} = - \ {\bf I} \ . 
\end{equation}
Note that the quantities $E_k$ and $E_z$ in \eqref{LE-all-states} and \eqref{LE-all-modes+}
are defined in the spaces of state variables and modal vectors, respectively. 
Therefore, 
in general, $E_k \ne E_z$.\\

\noindent {\bf Remark~1}. 
Without the additional normalization condition \eqref{ev-assumption-2}, the invariant definition of the Lyapunov modal energy takes the form
\[
E_{z_i} = |\boldsymbol{u}_i|^2 \int^{\infty}_0 z^*_i(t) \, z_i(t)\,dt \ .
\]

\noindent {\bf Remark~2}. 
The Lyapunov energies of states \eqref{Lener-k-state} and modes \eqref{Lener-i-mode+} depend on the dimensional units of the system variables. 
Using additional dimensional coefficients $\boldsymbol{c} = (c_1,\cdots, c_n)^T \in \mathbb{ R}^n$, these quantities can be represented in the form invariant 
with respect to a change of units 
\[
E_{x_k} = c^2_k \int^{\infty}_0 x^2_k(t)\,dt \  , \ \ \ 
E_{z_i} = |\boldsymbol{c}^T\boldsymbol{u}_i|^2 \int^{\infty}_0 z^*_i(t) \, z_i(t)\,dt \ .
\]
In the following, we assume that all state variables are suitably scaled so that $c_k = 1$ for all $k$.\\

On the basis of the definitions \eqref{Lener-k-state} and \eqref{Lener-i-mode+} 
of Lyapunov energies of the states and modes, we further introduce a novel concept of Lyapunov participation factors. \\

\subsection{Mode-in-state Lyapunov PF} 

The Lyapunov energy of each state $ E_ {x_k} $ in \eqref{Lener-k-state} 
can be partitioned into parts associated with individual eigenvalues 
using the decomposition \eqref{subgram+}-\eqref{decomp}
of Gramian ${\bf P}_{x_k}$ into the sub-Gramians 
\begin{equation} \label{Exk-decomp}
E_{x_k} = \sum^n_{i=1} E_{x_k i} \ , \ \text{where} \ \ E_{x_k i} = \boldsymbol{x}^T \widetilde{\bf P}_{x_k i} \, \boldsymbol{x} \ , \ 
{\bf P}_{x_k} = \sum^n_{i=1} \widetilde{\bf P}_{x_k i} \ ,
\end{equation}
and the sub-Gramians $\widetilde{\bf P}_{x_k i}$, which, according to \eqref{lyap-subgram}, satisfy the following Lyapunov equations
\[
 {\bf A}^* \, \widetilde{\bf P}_{x_k  i} + \widetilde{\bf P}_{x_k  i} \, {\bf A} = - \ \frac{1}{2}
 \Big( {\bf R}^*_i {\bf e}_k {\bf e}^T_k + {\bf e}_k {\bf e}^T_k {\bf R}_i \Big) \ .
\]
According to \eqref{subgram-eigen}-\eqref{pair-subgram-eigen}, for the sub-Gramians that we obtained 
\begin{gather} 
\widetilde{\bf P}_{x_k  i} = - \left\{(\boldsymbol{\varv}^T_i)^* (u^k_i)^* {\bf e}^T_k (\lambda^*_i {\bf I} + {\bf A})^{-1} \right\}_H \nonumber \\
= - \left\{ \sum^n_{j=1}\cfrac{(\boldsymbol{\varv}^T_i)^* (u^k_i)^* u^k_j \, \boldsymbol{\varv}^T_j}{\lambda_i^* +\lambda_j}\right\}_H \ . \label{Pxki}
\end {gather}
Therefore, Lyapunov energies of the states \eqref{Lener-k-state} can be
expressed through the eigenvectors and state variables as
\begin{gather}
E_{x_k} = {\boldsymbol{x}}^T {\bf P}_{x_k} {\boldsymbol{x}} =
 - {\boldsymbol{x}}^T \left\{ \sum^n_{i=1} ({\boldsymbol{\varv}}^T_i)^* (u^k_i)^*  
{\bf e}^T_k (\lambda^*_i {\bf I} +{\bf A})^{-1} \right\}_H {\boldsymbol{x}} \nonumber \\
=  - {\boldsymbol{x}}^T \left\{ \sum^n_{i,j=1} 
\frac{ ({\boldsymbol{\varv}}^T_i)^* (u^k_i)^* u^k_j \, {\boldsymbol{\varv}}^T_j} {\lambda^*_i + \lambda_j}
\right\}_H {\boldsymbol{x}} \ . \label{Pxk}
\end{gather}
This expression can also be obtained directly from the solution of the system \eqref{system0}
\[
x_k(t) = \sum^n_{i=1} u_i^k \, z_i(0) \,{\rm exp}\,(\lambda_i t) , \ \text{where} \ z_i(0) = \boldsymbol{\varv}^T_i \boldsymbol{x}_0, \ 
\boldsymbol{x}_0 = \boldsymbol{x}(0) .
\]
Substituting this in the definition \eqref{Lener-k-state}, we obtain the same result as in \eqref{Pxk}
\begin{gather}
E_{x_k}=\int^{\infty}_0 x^*_k(t) x_k(t) dt = \nonumber \\
= \int^{\infty}_0 
\Big( \sum\limits_{i=1}^n {\boldsymbol{x}}^T_0 ({\boldsymbol{\varv}}^T_i)^* (u^k_i)^* {\rm exp}\,(\lambda^*_i t) \Big) \cdot 
\Big( \sum\limits_{j=1}^n u^k_j \, {\boldsymbol{\varv}}^T_j {\boldsymbol{x}}_0 \,{\rm exp}\,(\lambda_j t) \Big) \, dt \nonumber \\
= {\boldsymbol{x}}^T_0 \left( \sum^n_{i,j=1} 
 ({\boldsymbol{\varv}}^T_i)^* (u^k_i)^* u^k_j \, {\boldsymbol{\varv}}^T_j \int^{\infty}_0 {\rm exp}\,((\lambda^*_i +\lambda_j) t) \, dt
\right) {\boldsymbol{x}}_0 \ .\nonumber 
\end{gather}

On the basis of \eqref{Exk-decomp}, we introduce the following definition. \\

\noindent {\bf Definition~1}.
{\it The mode-in-state Lyapunov participation factor (MISLPF)} is 
the relative participation of the $i$-th mode in the Lyapunov energy ${E}_{x_k}$ of the $k$-th state \eqref{Lener-k-state}, that is 
\begin{equation} \label{def-1}
e_{ki} =\cfrac{E_{x_k i}}{E_{x_k}} = \cfrac{\boldsymbol{x}_0^T \, \widetilde{\bf P}_{x_k i} \, \boldsymbol{x}_0}
{\boldsymbol{x}_0^T \,{\bf P}_{x_k} \boldsymbol{x}_0} \ \ .
\end{equation}

The MISLPFs defined in \eqref{def-1} clearly depend on 
the initial state vector $\boldsymbol{x}_0$ and the correlation among its components. 
Two different ways to take into account the initial conditions have been employed 
in the conventional definitions of participation factors (PFs). 
The original definition \eqref{PF-def} of the PFs and generalized participations (GPs) in \citep{perez:82} used specially selected initial conditions of the form
\begin{equation} \label{ini-cond-convent}
 \boldsymbol{x}_0 = {\bf e}_k \ \ \text{and} \ \  \boldsymbol{x}_0 = {\bf e}_l , \ l\ne k \ .
\end{equation}
In particular, for the first initial condition, we get 
PF $p_{ki}$, which $``$weights$"$ the participation of the $i$\nobreakdash-th mode and
initial $k$\nobreakdash-th state in the dynamics of the $k$\nobreakdash-th state. 
For the second initial condition, we get 
the GP $p_{kil}$, which $``$weights$"$ the participation of the $i$\nobreakdash-th mode
and initial $l$\nobreakdash-th state in the dynamics of the \mbox{$k$-th} state 
\begin{equation} \label{PF-dynamic} \nonumber
x_k(t) = \sum_ip_{ki} \, x_k(0) \,{\rm exp}\,(\lambda_i t) + \sum_{i,\,l\ne k}  p_{kil} \, x_l(0) \,{\rm exp}\,(\lambda_i t) \ . 
\end{equation}
Under the initial conditions \eqref{ini-cond-convent}, by analogy with the $p_{ki}$ and $p_{kil}$, the corresponding \
mode-in-state Lyapunov PFs and GPs can be calculated as
\begin{equation} \label{MISLPF-ini-cond-1}
{e}_{ki} = \dfrac{\big(\widetilde{\bf P}_{x_k i}\big)_{kk}}{ \sum^n_{i'=1} \big(\widetilde{\bf P}_{x_k i'}\big)_{kk}} \ , \ \
{e}_{kil}
= \dfrac{\big(\widetilde{\bf P}_{x_k i}\big)_{ll}}{ \sum^n_{i'=1} \big(\widetilde{\bf P}_{x_k i'}\big)_{ll}} \ \ .
\end{equation}
The difference from $p_{ki}$ and $p_{kil}$ is that the participations in \eqref{MISLPF-ini-cond-1} 
are taken not in relation to the state variable itself, but in relation to its Lyapunov energy. \\

Another approach was proposed in \citep{hash:09}, which 
was based on averaging over an uncertain set of initial conditions. The corresponding 
formula for calculation of the PFs was proposed 
for the spherically symmetric distribution of the initial conditions with respect to zero. 
Following this approach, from Definition~1, we obtain
the following formula for calculating the MISLPF:
\begin{equation}
{e}_{ki} = \cfrac{ \mathrm{trace} \, \big( \, \widetilde{\bf P}_{x_k i} \big) }{ \mathrm{trace} \, \big({\bf P}_{x_k}\big)} \ .
\end{equation}
Substituting here from \eqref{Pxki}, we obtain 
\begin{gather}
{e}_{ki} = 
\frac{ {\rm Re} \left\{ \mathrm{trace} \left( ({\boldsymbol{\varv}}^T_i)^* (u^k_i)^*  
{\bf e}^T_k (\lambda^*_i {\bf I} +{\bf A})^{-1} \right) \right\}}
{{\rm Re} \left\{ \mathrm{trace} \left( \sum^n_{i'=1} ({\boldsymbol{\varv}}^T_{i'})^* (u^k_{i'})^*  
{\bf e}^T_k (\lambda^*_{i'} {\bf I} +{\bf A})^{-1} \right) \right\}} \nonumber \\
= \frac{ {\rm Re} \left\{ 
 \sum\limits^n_{j=1} 
\cfrac{ {\boldsymbol{\varv}}^*_i \boldsymbol{\varv}_j \,(u^k_i)^* u^k_j } {\lambda^*_i + \lambda_j} \right\}}
{ {\rm Re} \left\{ \sum\limits^n_{i',j=1} 
\cfrac{ {\boldsymbol{\varv}}^*_{i'} \boldsymbol{\varv}_j (u^k_{i'})^* u^k_j} {\lambda^*_{i'}+ \lambda_j} \right\}} \ . \label{eki-evx}
\end{gather}

\subsection{State-in-mode Lyapunov PF}

The initial Lyapunov energy of the $i$-th mode \eqref{Lener-i-mode+}
can be expressed through the eigenvectors and initial state variables: 
 \begin{gather}
E_{z_i} = {\boldsymbol{x}}_0^T {\bf V}^* {\bf P}_{\boldsymbol{z} i} {\bf V} {\boldsymbol{x}}_0 =
\cfrac{{\boldsymbol{x}}_0^T ({\boldsymbol{\varv}}^T_i)^* \, {\boldsymbol{\varv}}^T_i {\boldsymbol{x}}_0}
{-2 \, {\rm Re} \{\lambda_i \} } , \ \text{where} \ 
{\bf P}_{\boldsymbol{z} i}=\cfrac{-{\bf e}_i \, {\bf e}^T_i}{2 \, {\rm Re} \{\lambda_i \} } \ .  \nonumber 
\end{gather}
This Lyapunov energy can be partitioned
into the parts corresponding to the state variables so that 
the contribution from each two state variables be divided in half between them.
Then, the contribution in $E_{z_i}$ from the $k$-th state is
\begin{gather} 
E_{z_i  k} =\frac{1}{2}  {\boldsymbol{x}}_0^T\Big( {\bf V}^* {\bf P}_{\boldsymbol{z} i} {\bf V} \, {\bf e}_k {\bf e}^T_k +
{\bf e}_k {\bf e}^T_k {\bf V}^* {\bf P}_{\boldsymbol{z} i} {\bf V} \Big) \,{\boldsymbol{x}}_0 \ . \nonumber \\
= \cfrac{{\boldsymbol{x}}_0^T ({\boldsymbol{\varv}}^T_i)^* \varv^k_i x^0_k + (\varv^k_i)^* x^0_k \, {\boldsymbol{\varv}}^T_i {\boldsymbol{x}}_0}
{-4 \, {\rm Re} \{\lambda_i \} } \ ,\label{Ezik}
\end{gather}
where $\varv^k_i$ and $x^0_k$ are the $k$-th components of ${\boldsymbol{\varv}}_i$ and ${\boldsymbol{x}}_0$, respectively. 
Introducing the notation 
\begin{gather} 
E_{z_i k} = {\boldsymbol{z}}_0^* \, \widetilde{\bf P}_{\boldsymbol{z} i k} \, \boldsymbol{z}_0 \ , 
\ \ \widetilde{\bf P}_{\boldsymbol{z} i k} = {\bf P}_{\boldsymbol{z} i} {\bf V} \, {\bf e}_k {\bf e}^T_k {\bf U}
 + {\bf U}^* {\bf e}_k {\bf e}^T_k {\bf V}^* {\bf P}_{\boldsymbol{z} i} \ , \nonumber 
\end{gather}
we can formulate the following definition. \\

\noindent {\bf Definition~2}.
{\it The state-in-mode Lyapunov participation factor (SIMLPF)} characterizes 
the relative participation of the $k$-th state in the Lyapunov energy ${E}_{z_i}$ of the $i$-th mode \eqref{Lener-i-mode+}, that is 
\begin{equation} \label{def-2}
\varepsilon_{ki} =\cfrac{E_{z_i k}}{E_{z_i}} =
\cfrac{{\boldsymbol{z}}_0^* \, \widetilde{\bf P}_{\boldsymbol{z} i k} \, \boldsymbol{z}_0}
{{\boldsymbol{z}}_0^*\, {\bf P}_{\boldsymbol{z} i} \, \boldsymbol{z}_0} 
= \cfrac{{\boldsymbol{x}}_0^T ({\boldsymbol{\varv}}^T_i)^* \varv^k_i x^0_k + (\varv^k_i)^* x^0_k \, {\boldsymbol{\varv}}^T_i {\boldsymbol{x}}_0}
{2 \, {\boldsymbol{x}}_0^T ({\boldsymbol{\varv}}^T_i)^* \, {\boldsymbol{\varv}}^T_i {\boldsymbol{x}}_0 } \ .
\end{equation}\\

This definition coincides with the modified definition for conventional SIMPF proposed in \citep{isk:19}. 
In this study, it was shown that Definition~2 in general
reproduces the results of the definition of SIMPFs proposed in \citep{hash:09}
for real eigenvalues, and rectifies the results in the case of complex eigenvalues. 
In particular, for the spherically symmetric distribution 
of the initial conditions with respect to zero, we obtain from Definition~2,
\begin{equation} \label{simpf+}
\varepsilon_{ki} = 
\cfrac{(\varv^k_i)^* \, \varv^k_i}{{{\left({\boldsymbol{\varv}}_i\right)}^* \, \boldsymbol{\varv}}_i } \ . 
\end{equation}
This coincides with \eqref{simpf} obtained in \citep{hash:09} 
for real eigenvalues and corrected in \citep{konoval:17} for the case of complex eigenvalues.
The obtained coincidence can be explained as follows.
The amplitude of the $i$-th mode and all its components 
(into whatever parts it is divided) change in time with the same exponential rate. 
Therefore, the ratio of the corresponding Lyapunov energies 
coincides with the ratio of the same $``$instantaneous$"$ energies at the initial moment. 
Thus, the state-in-mode LPFs coincide with the corresponding conventional PFs. 
In contrast, the mode-in-state LPFs in \eqref{def-1} are characterized by Lyapunov energies 
in the state space, which are determined by the exponential terms with different rates. 
Therefore, in contrast to \eqref{def-2}, there is a characteristic dependence in \eqref{eki-evx}
on pair combinations of the eigenvalues.\\

\noindent {\bf Remark}.
Although Definition~2 is invariant with respect to a change of units of the system state variables, the expression \eqref{simpf+} is not,
because it includes the assumption of spherical symmetry of the initial conditions. The expression \eqref{simpf+}, however, justifies its meaning 
when all state variables are independent and scaled so that they have the same variance.

\section{Modal analysis of Lyapunov energy of states}

\subsection{Modal contributions to Lyapunov energy of states}

The Lyapunov modal energy $E_z$ in \eqref{LE-all-modes+} is defined in the space of modal vectors~$\boldsymbol{z} \in \mathbb{C}^n$.
Nevertheless, the Lyapunov energy of states $E_x$ in \eqref{LE-all-states}
can also be considered from the point of view of modes, if we apply the transformation \eqref{z-vector}, that is
\begin{equation} \label{sumEi=sumEk+}
E_x = \boldsymbol{x}^T {\bf P}_{\boldsymbol{x}} \, \boldsymbol{x} =
\boldsymbol{x}^T {\bf V}^* {\bf U}^* {\bf P}_{\boldsymbol{x}} \, {\bf U} \, {\bf V} \, \boldsymbol{x} = 
\boldsymbol{z}^* \, \overline{\bf P}_{\boldsymbol{z}} \, \boldsymbol{z} = \overline{E}_z \ ,
\end{equation}
where 
$\overline{\bf P}_{\boldsymbol{z}} ={\bf U}^* {\bf P}_{\boldsymbol{x}} \, {\bf U}$, 
and  
$\boldsymbol{z} = {\bf V} \boldsymbol{x}$. 
Substituting 
the representation \eqref{A=ULV} of the matrix ${\bf A}$ into \eqref{Lyap-states}
and multiplying the resulting equation from the right by the matrix ${\bf U}$ 
 and from the left by the matrix ${\bf U}^*$, we obtain that $\overline{\bf P}_{\boldsymbol{z}}$ satisfies 
 the following Lyapunov equation:  
\begin{equation} \label{Lyap-eq-modes}
{\bf \Lambda}^* \, \overline{\bf P}_{\boldsymbol{z}} +
\overline{\bf P}_{\boldsymbol{z}} \, {\bf \Lambda} = - \, {\bf U}^*  {\bf U} \ . \ 
\end{equation}

The joint contribution to $E_x$ from $i$th and $j$th modes can be divided equally between these modes. 
Then, a contribution of the $i$-th mode to $E_x$ can be defined as
\begin{equation} \label{Ezi-contrib}
\overline{E}_{z_i} =\frac{1}{2} \sum^n_{j=1} \Big(z^*_i \, (\overline{\bf P}_{\boldsymbol{z}})_{ij} \, z_j + z^*_j \, (\overline{\bf P}_{\boldsymbol{z}})_{ji} \, z_i \Big) \ ,
\end{equation}
where $(\overline{\bf P}_{\boldsymbol{z}})_{ij}$ and $(\overline{\bf P}_{\boldsymbol{z}})_{ji}$ 
are the elements of the matrix $\overline{\bf P}_{\boldsymbol{z}}$. Substituting their values from \eqref{Lyap-eq-modes}, we obtain
\begin{equation} \label{Ei-expression}
\overline{E}_{z_i} = - \, {\rm Re} \left\{ \sum^n_{j=1} \frac{\boldsymbol{u}^*_i \, \boldsymbol{u}_j}{\lambda_i^* +\lambda_j} z^*_i \, z_j \right\} \ .
\end{equation}

On the other hand, one can apply to \eqref{Lyap-eq-modes} the decomposition
(\ref{subgram+})-(\ref{decomp}) of $\overline{\bf P}_{\boldsymbol{z}}$ into the sub-Gramians:  
\begin{gather} 
\overline{\bf P}_{\boldsymbol{z}} = \sum^n_{i=1} \overline{\bf P}_{\boldsymbol{z} i} \ , \ \ \text{where} \nonumber \\
\overline{\bf P}_{\boldsymbol{z}i} = - \left\{{\bf e}_i {\bf e}^T_i {\bf U}^*  {\bf U}  (\lambda^*_i {\bf I} + {\bf \Lambda})^{-1} \right\}_H 
= -  \left\{ \sum^n_{j=1} \frac{{\bf e}_i \, \boldsymbol{u}^*_i  \boldsymbol{u}_j \,{\bf e}_j^T}{\lambda_i^* +\lambda_j} \right\}_H . \label{Pzi}
\end{gather}
Therefore, comparing this with \eqref{Ei-expression}, we can define the {\it
modal contribution (MC) of the $i$-th mode to the Lyapunov energy of states} 
using the sub-Gramian $\overline{\bf P}_{\boldsymbol{z} i}$ as 
\begin{gather}
\overline{E}_{z_i} = {\boldsymbol{z}}^* \overline{\bf P}_{\boldsymbol{z} i} \, \boldsymbol{z} \ , \ \ \text{where} 
\ \ \ \ \boldsymbol{z} = {\bf V} \boldsymbol{x} \ ,\nonumber \\
{\bf \Lambda}^* \, \overline{\bf P}_{\boldsymbol{z} i} +
\overline{\bf P}_{\boldsymbol{z} i} \, {\bf \Lambda} = - \frac{1}{2}
\Big( {\bf U}^*{\bf U} \, {\bf e}_i {\bf e}^T_i + {\bf e}_i {\bf e}^T_i  {\bf U}^*{\bf U} \Big) \ .  \label{Lener-i-mode}
\end{gather}

Values $E_{z_i}$ in \eqref{Lener-i-mode+} characterize modal energies in the space of modal vectors. 
By definition, they are always positive. 
In contrast, the quantities $\overline{E}_{z_i}$ in \eqref{Lener-i-mode} determine the contributions 
of the modes to the Lyapunov energy accumulated in state variables. 
As can be seen from \eqref{Ei-expression}, they can be negative if there are modes 
with sufficiently large amplitudes and negative correlation with a given mode.
It follows from \eqref{Ezi-contrib} that 
\begin{equation}
\overline{\bf P}_{\boldsymbol{z} i} = \frac{1}{2} \Big(\overline{\bf P}_{\boldsymbol{z}} \, {\bf e}_i {\bf e}^T_i 
+ {\bf e}_i {\bf e}^T_i \, \overline{\bf P}_{\boldsymbol{z}} \Big) \ .
\end{equation}
Because the matrix $\overline{\bf P}_{\boldsymbol{z} i}$ is self-adjoint,  $\big(\overline{E}_{z_i}\big)^* = \overline{E}_{z_i} $, 
i.e., the MCs 
are always real.
In addition, the MCs 
of the modes corresponding to a pair of complex conjugate eigenvalues are the same.\\

\noindent {\bf Proposition~1}.
Let $i$ and $i*$ be the indexes  
associated with eigenvalues $\lambda_i$  and $\lambda^*_i$, respectively. 
Then, the corresponding modal contributions 
are the same:  
\begin{equation} \nonumber
\overline{E}_{z_{i*}} = \overline{E}_{z_i} \ .
\end{equation}
\noindent {\bf Proof.}
For arbitrary $j$, we consider in the expression \eqref{Lener-i-mode} for $\overline{E}_{z_i}$ a term 
\begin{equation} \label{prop-1-1}
z_i^* \, \overline{\bf P}_{\boldsymbol{z} i} \, z_j + z_j^*\, \overline{\bf P}_{\boldsymbol{z} i}\,  z_i = -\frac{1}{2} \left(
 \cfrac{z_i^* \boldsymbol{u}^*_i  \boldsymbol{u}_j z_j }{\lambda_i^* +\lambda_j} 
+\cfrac{z_j^* \boldsymbol{u}^*_j  \boldsymbol{u}_i z_i }{\lambda_j^* +\lambda_i} \right)_H \ . 
\end{equation}
Choose $j*$ such that if $j$ corresponds to a real eigenvalue, then $j*=j$, 
and if $j$ corresponds to a complex eigenvalue $\lambda_j$, then $j*$ corresponds to a complex eigenvalue $\lambda^*_j$.
Then, one can verify that in the expression for $\overline{E}_{z_{i *}}$, there is the term
\begin{equation} \label{prop-1-2}
z_{i*}^*\, \overline{\bf P}_{\boldsymbol{z} i*} \,z_{j*} + z_{j*}^*\, \overline{\bf P}_{\boldsymbol{z} i*} \, z_{i*}  \ ,
\end{equation}
which equals to \eqref{prop-1-1} under the assumption \eqref{ev-assumption}.
By the one-to-one correspondence of terms \eqref{prop-1-1} and \eqref{prop-1-2}, we obtain,
${\boldsymbol{z}}^* \,\overline{\bf P}_{\boldsymbol{z} i} \, \boldsymbol{z} =
 {\boldsymbol{z}}^* \,\overline{\bf P}_{\boldsymbol{z} i*} \, \boldsymbol{z}$. $\square$ \\

The total contribution of all modes is equal to the Lyapunov energy of all states: 
\begin{equation} \label{LE-all-modes}
\sum^n_{i=1} \overline{E}_{z_i} = {\boldsymbol{z}}^* \overline{\bf P}_{\boldsymbol{z}} \, {\boldsymbol{z}} = E_x \ , \ \text{where} \ \ 
\overline{\bf P}_{\boldsymbol{z}} = \sum^n_{i=1} \overline{\bf P}_{\boldsymbol{z} i} 
\end{equation}
is the solution of the Lyapunov equation \eqref{Lyap-eq-modes}.
From \eqref{Pzi} and \eqref{Pxki}, one can also verify that the modal contribution of the $i$-th mode is 
the sum of its contributions to all states  
\[
\overline{E}_{z_i} = {\boldsymbol{x}}^T {\bf V}^* \overline{\bf P}_{z_i} {\bf V} \boldsymbol{x} =
\sum^n_{k=1}  {\boldsymbol{x}}^T \widetilde{\bf P}_{x_k  i} \boldsymbol{x} = \sum^n_{k=1} E_{x_k i} \ .
\]

\subsection{Lyapunov modal interactions}

Unlike state variables or signals, their Lyapunov energies can be partitioned 
not only into parts corresponding to individual eigenvalues, 
but also into parts corresponding to pair combinations of eigenvalues, 
i.e., the solution of Lyapunov equation \eqref{lyapunov} can be decomposed by 
\eqref{subgram+}-\eqref{decomp}:
\begin{equation} \nonumber 
{\bf P} = \sum_{i,j=1}^n {\bf P}_{ij} = - \sum\limits_{i,j=1}^n \left\{ \frac{{\bf R}^*_i \ {\bf Q} \ {\bf R}_j}{\lambda^*_i+\lambda_j} \right\}_H \ .
\end{equation}
This allows us to characterize 
the pairwise modal interactions in the system 
by the Lyapunov energies produced in the system states by the corresponding pairwise mode combinations.\\

\noindent {\bf Definition~3}.
{\it The Lyapunov modal interaction energy (LMIE)} of the $i$-th and $j$-th modes in the system \eqref{system0} is
\begin{gather}
\overline{E}_{z \,i j} = {\boldsymbol{x}}^T {\bf P}_{(i j)} \,\boldsymbol{x} \ , \ \text{where} \  \
{\bf P}_{(i j)} = - 
\frac{1}{2} \left\{ \cfrac{{\bf R}^*_i \, {\bf R}_j}{\lambda^*_i + \lambda_j } + \cfrac{{\bf R}^T_i \, {\bf R}_j}{\lambda_i + \lambda_j }  \right\}_{H} 
\nonumber \\
= - 
\frac{1}{2} \left\{ \frac{ ({\boldsymbol{\varv}}^T_i)^* \boldsymbol{u}^*_i \, \boldsymbol{u}_j \, {\boldsymbol{\varv}}^T_j} {\lambda^*_i + \lambda_j} 
+ \frac{ {\boldsymbol{\varv}}_i \ \boldsymbol{u}^T_i \, \boldsymbol{u}_j \, {\boldsymbol{\varv}}^T_j} {\lambda_i + \lambda_j} 
\right\}_H 
 \ . \label{def-3}
\end{gather}\\

Because the matrix ${\bf P}_{(i j)}$ is self-adjoint, the LMIEs are always real, 
i.e., $\overline{E}_{z \,i j} = {\boldsymbol{x}}^T {\rm Re} \left\{{\bf P}_{(i j)} \right\} \boldsymbol{x}$.
In addition, for any matrix ${\bf C}$, the equality ${\rm Re}\left\{\left\{{\bf C}\right\}_H\right\} = {\rm Re}\left\{\left\{{\bf C}\right\}_{SYM}\right\}$ holds, 
where $\{{\bf C}\}_{SYM} = \frac{1}{2}( {\bf C} +{\bf C}^T)$.
Therefore, the definition of LMIE implies that
\begin{gather}
\overline{E}_{z j\, i} = - \frac{1}{2} {\boldsymbol{x}}^T {\rm Re}\left\{ \left\{ \cfrac{{\bf R}^*_j \, {\bf R}_i }{\lambda^*_j + \lambda_i } \right\}_H  
+ \left\{\cfrac{{\bf R}^T_j \, {\bf R}_i}{\lambda_j + \lambda_i }  \right\}_{SYM} \right\} \boldsymbol{x} \nonumber \\
= - \frac{1}{2} {\boldsymbol{x}}^T {\rm Re}\left\{ \left\{ \cfrac{ {\bf R}^*_i \, {\bf R}_j}{\lambda^*_i + \lambda_j } \right\}_H 
+ \left\{\cfrac{{\bf R}^T_i \, {\bf R}_j}{\lambda_i + \lambda_j }  \right\}_{SYM} \right\} \boldsymbol{x} =
\overline{E}_{z\,i j} \ . \nonumber
\end{gather}
Furthermore, the symmetrization in \eqref{def-3} ensures 
that the LMIEs associated with the complex conjugate eigenvalues $\lambda_i$ and $\lambda^*_i$ are always the same, that is
\begin{equation}
\overline{E}_{z\,i* j} = \overline{E}_{z \,i j} \ , 
\end{equation}
where $i$ and $i*$ are the indexes 
associated with eigenvalues $\lambda_i$ and $\lambda^*_i$, respectively. 
The matrix ${\bf P}_{(i j)}$ in \eqref{def-3} is the solution of the Lyapunov equation
\begin{equation} \nonumber
  {\bf A}^* \, {\bf P}_{(i j)} + {\bf P}_{(i j)} \, {\bf A} = - \ \frac{1}{4}\,\Big({\bf R}^*_i {\bf R}_j
  +{\bf R}^T_i {\bf R}_j + {\bf R}^*_j {\bf R}_i + {\bf R}^*_j \big({\bf R}^T_i\big)^* \Big)  . 
\end{equation}
It follows from \eqref{Ei-expression}, \eqref{Pzi}, and \eqref{def-3} that 
\begin{equation}
\frac{1}{2}\,\Big(\overline{E}_{z_i}+\overline{E}_{z_{i*}}\Big) = \sum^n_{j=1} \overline{E}_{z\,i j} \ , \ \ 
\frac{1}{2}\,{\bf V}^* \Big( \overline{\bf P}_{\boldsymbol{z} i} + \overline{\bf P}_{\boldsymbol{z} i*}\Big){\bf V} = \sum^n_{j=1} {\bf P}_{(i j)} \ . \nonumber
\end{equation}
In view of Proposition~1, this implies that the MC of the mode is 
a sum of the Lyapunov energies 
corresponding to its modal interactions with other modes, that is  
\begin{equation}
\overline{E}_{z_i} = \sum^n_{j=1} \overline{E}_{z \,i j} \ , \ \ 
{\bf V}^* \overline{\bf P}_{\boldsymbol{z} i} {\bf V} = \sum^n_{j=1} {\bf P}_{(i j)} \ . \nonumber
\end{equation}

When analyzing small-signal stability in a linearized system, it would be convenient to have general indicators 
that would characterize the participation of the modes and their interaction 
in terms of accumulated Lyapunov energy, but would not depend on state variables, in which the initial perturbation was made. 
For this purpose, we consider a probabilistic description of the uncertainty in the initial condition
and assume that the components of the initial condition vector $x^0_1, x^0_2 ,\ldots, x^0_n$ 
are distributed independently with zero mean and unit variance, that is
\begin{equation} \label{assume-spheric-symmetry}
{\rm E} \, \big\{ x^0_k x^0_{k'}\big\} = \delta_{kk'} \ ,
\end{equation}
where the expectation operator ${\rm E}\left\{ \cdot \right\}$ is evaluated using some assumed joint probability function $f(\boldsymbol{x}_0)$ for the initial condition uncertainty, and $\delta_{kk'}$ is the Kronecker delta.\\

\noindent {\bf Proposition~2}.
Under assumption \eqref{assume-spheric-symmetry},
the averaged Lyapunov MCs 
$\overline{E}_{z_i}$ and modal interaction energies $\overline{E}_{z\,i j}$ can be computed as
 \begin{gather} 
{\rm E} \, \Big\{ \overline{E}_{z_i} \Big\} = - \ {\rm Re} \left\{ \mathrm{trace} \left(  {\bf R}^*_i \, (\lambda^*_i {\bf I} + {\bf A})^{-1}  \right) \right\} \ , \nonumber \\
{\rm E} \, \Big\{ \overline{E}_{z \,i j} \Big\}= - \ \frac{1}{2} \, {\rm Re} \left\{ \frac{  \mathrm{trace} \left( {\bf R}^*_i  {\bf R}_j \right) }{\lambda^*_i+\lambda_j} 
+  \frac{  \mathrm{trace} \left( {\bf R}^T_i {\bf R}_j \right) }{\lambda_i+\lambda_j} \right\}  \ . \label{Ezij-averaged}
 \end{gather}
 
\noindent {\bf Proof.}
Under the assumption \eqref{assume-spheric-symmetry}, the first expression follows from  \eqref{Pzi}, \eqref{Lener-i-mode}, 
and the formula ${\bf R}_i = \boldsymbol{u}_i \boldsymbol{\varv}^T_i$ for the residue 
of the resolvent of the matrix ${\bf A}$. The second expression follows directly from the definition \eqref {def-3}. 
We also take into account the equality ${\rm trace}\,\{{\bf C}\}_H = {\rm Re}\{{\rm trace} \,{\bf C}\}$. $\square$ \\

Based on the result of Proposition~2, the following indicator can be proposed as a convenient measure of the relative participation of the modes in each other.\\

\noindent {\bf Definition~4}.
{\it The Lyapunov modal interaction factor (LMIF)} of the $j$-th mode in the $i$-th mode is 
\begin{equation} \label{def-4}
LMIF_{ij} = \cfrac{{\rm E} \, \Big\{ \overline{E}_{z \,i j} \Big\}} {\sum^n_{j'=1} \left|{\rm E} \, \Big\{ \overline{E}_{z \,i j'} \Big\}\right| } \ \ ,
\end{equation}
where ${\rm E} \, \Big\{ \overline{E}_{z \,i j} \Big\}$ is defined in \eqref{Ezij-averaged}.\\

These indicators make practical sense if the state variables can be scaled in such a way 
that their disturbances have approximately the same importance 
as the point of view of analyzing small signals in the system.
This somewhat specific interpretation is compensated, however, by the fact that the LMIFs
 do not have an explicit dependence on individual state variables.
 
\subsection{Lyapunov pair PFs}
 
In order to analyze modal interactions in connection with specific state variables, 
it is necessary to introduce Lyapunov energies and participation factors 
that correspond not only to individual eigenvalues but also to their pair combinations. 
According to \eqref{subgram+}, \eqref{decomp}, \eqref{lyap-subgram}, and \eqref{Lener-k-state}, 
{\it Lyapunov energy of the $k$-th state variable associated with a pair of $i$-th and $j$-th modes}
can be defined as
\begin{gather}
E_{x_k i j} = \boldsymbol{x}^T {\bf P}_{x_k i j} \, \boldsymbol{x} \ , \ \ \
{\bf P}_{x_k} = \sum^n_{i,j=1} {\bf P}_{x_k i j} \ , \ \ \ \text{where} \nonumber \\
{\bf P}_{x_k i j} = - \frac{1}{2}\left\{ \cfrac{{\bf R}^*_i {\bf e}_k {\bf e}^T_k {\bf R}_j}{\lambda_i^* +\lambda_j}
+ \cfrac{{\bf R}^T_i {\bf e}_k {\bf e}^T_k {\bf R}_j}{\lambda_i +\lambda_j}\right\}_H \nonumber \\
= -  \frac{1}{2} \left\{ \cfrac{(\boldsymbol{\varv}^T_i)^* (u^k_i)^* u^k_j \, \boldsymbol{\varv}^T_j}{\lambda_i^* +\lambda_j}
+ \cfrac{\boldsymbol{\varv}_i \, u^k_i \, u^k_j \, \boldsymbol{\varv}^T_j}{\lambda_i +\lambda_j} \right\}_H \ , \label{Pxkij}
\end{gather} 
and the symmetrized sub-Gramians ${\bf P}_{x_k i j}$ satisfy the Lyapunov equations
\begin{gather}
 {\bf A}^* \, {\bf P}_{x_k i j} + {\bf P}_{x_k i j} \, {\bf A} = \nonumber \\
 - \ \frac{1}{4}\,\Big({\bf R}^*_i {\bf e}_k {\bf e}^T_k {\bf R}_j + {\bf R}^*_j {\bf e}_k {\bf e}^T_k {\bf R}_i
+ {\bf R}^T_i {\bf e}_k {\bf e}^T_k {\bf R}_j + {\bf R}^*_j {\bf e}_k {\bf e}^T_k ({\bf R}^T_i)^* \Big)  \ . \nonumber
\end{gather}
By this definition, we have $E_{x_k i j} = E_{x_k j i} $ and $E_{x_k i^* j} = E_{x_k i j} $. 
Moreover, from \eqref{def-3} and \eqref{Pxkij} it follows that
\begin{equation} \label{Ezij=sumExkij}
\sum^n_{k=1} E_{x_k i j}  = \overline{E}_{z \,i j} \ . 
\end{equation}

Similar to \eqref{Pxkij},  
Lyapunov energy of the $i$-th mode associated with a pair of $k$-th and $l$-th states can be defined as
\begin{gather}
E_{z_{i} kl} = {\boldsymbol{z}}^* {\bf P}_{z_i kl} \,{\boldsymbol{z}}  \nonumber \\
= \frac{1}{2}  {\boldsymbol{x}}^T\Big( {\bf e}_l {\bf e}^T_l {\bf V}^* {\bf P}_{\boldsymbol{z} i} {\bf V} \, {\bf e}_k {\bf e}^T_k +
{\bf e}_k {\bf e}^T_k {\bf V}^* {\bf P}_{\boldsymbol{z} i} {\bf V} \, {\bf e}_l {\bf e}^T_l \Big) \,{\boldsymbol{x}}  \nonumber \\
= \cfrac{(\varv^l_i)^* x_l \, \varv^k_i x_k + (\varv^k_i)^* x_k \, \varv^l_i x_l }{-4 \, {\rm Re} \{\lambda_i \} } \ . \label{Eikl}
\end{gather} 

On the basis of Lyapunov energies \eqref{Pxkij} and \eqref{Eikl}, we can introduce the corresponding Lyapunov PFs.\\

\noindent {\bf Definition~5}.
{\it The pair MISLPF} is 
the relative pairwise participation of $i$-th and $j$-th modes 
in the Lyapunov energy \eqref{Lener-k-state} of the $k$-th state
\begin{equation} \label{def-5-1}
\tilde{e}_{k(ij)} = \cfrac{E_{x_k i j}}{E_{x_k}} =  \cfrac{\boldsymbol{x}^T {\bf P}_{x_k i j} \, \boldsymbol{x}}
{\boldsymbol{x}^T {\bf P}_{x_k}  \boldsymbol{x}} \ \ .
\end{equation}
{\it The pair SIMLPF} is  
the relative participation of $k$-th and $l$-th states in the
Lyapunov energy \eqref{Lener-i-mode+} of the $i$-th mode
\begin{equation} \label{def-5-2}
\tilde{\varepsilon}_{i(kl)} = \cfrac{E_{z_{i} kl}}{E_{z_i}} =  \cfrac{{\boldsymbol{z}}^* {\bf P}_{z_i kl} \,{\boldsymbol{z}}}
{{\boldsymbol{z}}^* {\bf P}_{z_i} {\boldsymbol{z}}} \ \ .
\end{equation} \\

The first indicator \eqref{def-5-1} shows 
the state variables having Lyapunov energies that are most sensitive to a particular modal interaction.
The second indicator \eqref{def-5-2} shows the pairs of state variables that produce the particular mode. 
Definition~5 is consistent with the previous definitions in the sense that the following relations hold between the Lyapunov PFs: 
\begin{equation}
\sum^n_{j=1} \tilde{e}_{k(ij)} = e_{ki}  , \ 
\sum^n_{l=1} \tilde{\varepsilon}_{i(kl)} = \varepsilon_{ki} , \ 
\sum^n_{i,j=1} \tilde{e}_{k(ij)} = \sum^n_{k,l=1} \tilde{\varepsilon}_{i(kl)}= 1 \ . \nonumber 
\end{equation}

Under the conventional initial conditions \eqref{ini-cond-convent}, 
the corresponding pair mode-in-state Lyapunov PFs and GPs can be calculated similar 
to \eqref{MISLPF-ini-cond-1} as
\begin{equation} \label{pair-MISLPF-ini-cond-1}
\tilde{e}_{k(ij)} = \dfrac{\big({\bf P}_{x_k i j}\big)_{kk}}{ \big({\bf P}_{x_k }\big)_{kk}} \ , \ \
\tilde{e}_{k(ij)l}
= \dfrac{\big( {\bf P}_{x_k i j}\big)_{ll}}{\big({\bf P}_{x_k}\big)_{ll} } \ \ .
\end{equation}

The calculation formulas for pair Lyapunov PFs can also be easily obtained 
for the spherically symmetric distribution of the initial conditions with respect to zero 
using \eqref{Pxkij}, \eqref{Eikl}, \eqref{Lener-i-mode+}, and \eqref{Pxk}:
\begin{gather} 
\tilde{e}_{k(ij)} = \cfrac{{\rm trace} \,(\,{\bf P}_{x_k i j})}{{\rm trace} \, (\,{\bf P}_{x_k})} 
= \frac{ \cfrac{1}{2} \, {\rm Re} \left\{  
\cfrac{ {\boldsymbol{\varv}}^*_i {\boldsymbol{\varv}}_j  \, (u^k_i)^* u^k_j }{\lambda^*_i + \lambda_j} + 
\cfrac{{\boldsymbol{\varv}}_i \, {\boldsymbol{\varv}}_j  \, u^k_i \, u^k_j }{\lambda_i + \lambda_j}  \right\}}
{ {\rm Re} \left\{  \sum\limits^n_{i',j'=1} 
\cfrac{ {\boldsymbol{\varv}}^*_{i'} {\boldsymbol{\varv}}_{j'}  \, (u^k_{i'})^* u^k_{j'}} {\lambda^*_{i'}+ \lambda_{j'}} \right\}} , \label{ekij-sym-ini} \\
\tilde{\varepsilon}_{i(kl)} = \delta_{kl} \varepsilon_{ki}\ . \ \nonumber
\end{gather}

Although the pair SIMLPF \eqref{def-5-2} is in some sense dual to the pair MISPF \eqref{def-5-1}, it 
does not seem to be a meaningful indicator.
As a more conceptual indicator, we can consider the state
participation in LMIEs, which characterizes 
the relative participation of the $k$-th state in the LMIE \eqref{def-3} of the $i$-th and $j$-th modes as
\begin{gather}
\overline{e}_{k(ij)}=
\cfrac{\overline{E}_{z_{i j} k}} {\overline{E}_{z \,i j} } = 
\cfrac{{\boldsymbol{x}}^T {\bf P}_{(i j)k} \,{\boldsymbol{x}}}
{\boldsymbol{x}^T {\bf P}_{(i j)} \,\boldsymbol{x}}
 \ , \ \ \text{where} \label{LMIE-pf} \\
\overline{E}_{z_{ij} k} = {\boldsymbol{x}}^T {\bf P}_{(i j)k} \,{\boldsymbol{x}} \ , \ 
{\bf P}_{(i j)} = \sum^n_{k=1} {\bf P}_{(i j) k} \ , \ \overline{E}_{z\,i j} = \sum^n_{k=1} \overline{E}_{z_{ij} k} \ , \nonumber \\
{\bf P}_{(i j) k} = \frac{1}{2} \, \Big( {\bf P}_{(i j)} \, {\bf e}_k {\bf e}^T_k + {\bf e}_k {\bf e}^T_k {\bf P}_{(i j)} \Big) \ . \nonumber
\end{gather}
This indicator shows the state variables that produce a given modal interaction.
Indicator \eqref{def-5-1} can be useful for selecting state variables that are most sensitive for identifying a specific interaction, 
while indicator \eqref{LMIE-pf} can be used to select the state variables for damping this interaction. 

\section{Lyapunov PFs in the selective modal analysis}

In this section, we present some characteristic properties of Lyapunov PFs that 
highlight the potential advantages of Lyapunov modal analysis.

\subsection{Relation between Lyapunov PFs and conventional PFs}

We compare conventional PFs \eqref{PF-def} and MISLPFs of Definition~1 
under conventional initial conditions \eqref{ini-cond-convent}. In this case, the conventional coefficients $p_{ki}$ and $p_{kil}$ 
correspond to the Lyapunov PFs from \eqref{MISLPF-ini-cond-1} and pair Lyapunov PFs from \eqref{pair-MISLPF-ini-cond-1}.\\

\noindent {\bf Property~1}.
{\it Under the conventional initial conditions} \eqref{ini-cond-convent}, {\it MISLPFs in Definition}~1 \eqref{def-1}
{\it and pair MISLPFs in Definition}~5 \eqref{def-5-1}
{\it can be represented through the conventional PFs} \eqref{PF-def} {\it and corresponding eigenvalues as follows: }
\begin{gather}  
e_{ki} = -  \frac{1}{{\bf e}^T_k {\bf P}_{x_k} {\bf e}_k} \left\{ p^*_{ki} {\bf e}^T_k (\lambda^*_i {\bf I} +{\bf A})^{-1} {\bf e}_k \right\}_H \ , \nonumber \\
\tilde{e}_{k(ij)} = - \frac{1}{2 \, {\bf e}^T_k {\bf P}_{x_k} {\bf e}_k}
\left\{ \frac{p^*_{ki} p_{kj}}{\lambda^*_i+\lambda_j} + \frac{p_{ki} p_{kj}}{\lambda_i+\lambda_j} \right\}_H 
. \label{prop-1} 
\end{gather}

\noindent {\bf Proof.}
Substitute \eqref{Pxki} and \eqref{Pxkij} into the definitions \eqref{MISLPF-ini-cond-1} and \eqref{pair-MISLPF-ini-cond-1}, 
for the MISLPF and pair MISLPF, respectively, obtained under the conventional initial conditions \eqref{ini-cond-convent}. 
Then, we take into account the definition of the conventional PFs \eqref{PF-def}. $\square$ \\

While the conventional PFs characterize the participation of modes in {\it the
instantaneous dynamics} of the state variables, the MISLPFs and pair MISLPFs characterize their participation in {\it the integrated energies} accumulated in the state variables.
Thus, in contrast to the conventional PFs, new indicators take into account the time responses of
specific devices, which are reflected by the dependence on the eigenvalues in the denominators 
in \eqref{prop-1}. \\

\subsection{LPFs indicate the distance from the stability boundary}

Suppose that the matrix ${\bf A}$ smoothly changes depending on the parameter $\gamma$ 
so that the following assumption is satisfied. \\

\noindent {\bf Assumption~1}. All conventional PFs \eqref{PF-def} remain limited, i.e.
\begin{equation} 
\forall \ k, i \ \ \ \exists \ C \in {\mathbb R}: \ \ \ |p_{ki} | < C \ . 
\end{equation}

Since the conventional PFs represent the sensitivities of eigenvalues, this assumption
means that the spectrum of the dynamic matrix changes smoothly with respect to $\gamma$. Then,
the characteristic property of LPFs is established by the following proposition.\\

\noindent {\bf Property~2}.
{\it Let the matrix ${\bf A}$ of the stable system} (\ref{system}) {\it change under Assumption}~1 {\it in such a way that
the $i$-th eigenvalue tends to the imaginary axis from the left, while other 
eigenvalues are limited: }
\begin{equation}
{\rm Re} \, (\lambda_i) \to -0 \ , \ \ {\rm Re}\, (\lambda_j) < - \alpha < 0, \ \ \lambda_j\ne \lambda_i , \lambda^*_i .
\end{equation} 
{\it Then, under conventional initial conditions \eqref{ini-cond-convent}, $E_{z i} \to + \infty$, 
and there is at least one state variable $k$ such that the following limiting relations are satisfied for MISLPFs:}
\begin{gather} \nonumber
\text{if} \ \lambda_i \ \text{is real}: \ E_{{x_k}i} \to +\infty , \ \ e_{ki} \to 1 \ , \ \ e_{kj} \to 0 \ , \ \ j\ne i \ ; \nonumber \\
\text{if} \  \lambda_i \ \text{is complex}: \ E_{{x_k}i}, E_{{x_k}i^*} \to +\infty , \\
\ e_{ki}, e_{ki^*} \to 0.5 \ , \ \ e_{kj} \to 0 \ , \ \ j\ne i, \, i^* \ . \nonumber
\end{gather}

\noindent {\bf Proof.} 
Consider the case when $\lambda_i$ is real.
Choose a state variable $k$ such that $|p_{ki}| \ne 0$. 
Note that under conventional initial conditions \eqref{ini-cond-convent}, $E_{x_k} = {\bf e}^T_k {\bf P}_{x_k} {\bf e}_k$.
Then, it follows from \eqref{prop-1} that 
\begin{gather}
E_{x_k} \tilde{e}_{kii} = - \left\{ \frac{p^*_{ki} p_{ki}}{\lambda^*_i+\lambda_i} \right\}_H = \frac{|p_{ki}|^2}{-2Re(\lambda_i)} 
\to +\infty \ , \nonumber \\
\forall j\ne i: \ |E_{x_k} \tilde{e}_{kij}| 
< \frac {|p^*_{ki}| \, |p_{kj}|}{-Re(\lambda_i) - Re (\lambda_{j})} < \frac{C^2}{\alpha} \ , \nonumber \\
\forall j,j' \ne i: \ |E_{x_k} \tilde{e}_{kjj'}| 
< \frac {|p^*_{kj}| \, |p_{kj'}|}{-Re(\lambda_{j'}) - Re (\lambda_j)} < \frac{C^2}{2\alpha} \ . \nonumber 
\end{gather}
From here, we obtain 
\begin{gather}
E_{{x_k}i} = E_{x_k} e_{ki} = E_{x_k} \sum_{j'=1}^n \tilde{e}_{kij'} \to +\infty \ , \nonumber \\
\forall j\ne i: \ E_{{x_k}j} = E_{x_k} e_{kj} = E_{x_k} \sum_{j'=1}^n \tilde{e}_{kjj'} <  \frac{(n+1)\, C^2}{2\alpha} \ , \nonumber \\
e_{ki} = \frac{E_{{x_k}i}}{\sum_{j'=1}^n E_{{x_k}j'}} \to 1 , \ \ \forall j\ne i: e_{kj} = \frac{E_{{x_k}j}}{\sum_{j'=1}^n E_{{x_k}j'}} \to 0 . \nonumber 
\end{gather}
The case of complex $\lambda_i$ is treated similarly.
$\square$ \\

It follows from Property~2 that, in contrast to conventional PFs, the values of LPFs  
depend on the distance of the corresponding eigenvalues from the stability boundary.  
Therefore, the proposed indicators account for the risk of stability loss associated with
individual weakly stable modes, i.e. modes having low decay rates. 
In expression \eqref{prop-1}, this is reflected by the
dependence of the denominator on $\lambda_i$. The closer the $i$-mode is to
the stability boundary, the smaller the value of the real part of $\lambda_i$ is, and
the greater the values of the corresponding LPFs and Lyapunov energies are. 

\subsection{Pair LPFs indicate resonant modal interactions}

Consider a stable system (\ref{system}) with two close eigenfrequencies $\omega_1 \approx \omega_2$
corresponding to the eigenvalues $\lambda_{1, 1^*} = - \alpha_1 \pm i \omega_1$ and $\lambda_{2, 2^*} = - \alpha_2 \pm i \omega_2$.
Suppose that the matrix ${\bf A}$ smoothly changes depending on the parameter $\gamma$ 
so that the following assumptions are satisfied. \\

\noindent {\bf Assumption~2}. The damping for $\lambda_{1}, \lambda_2$ is small, i.e. 
 \begin{equation} \label{assump-2}
\alpha_1 + \alpha_2 << \omega_1 + \omega_2 \ .
\end{equation}

\noindent {\bf Remark.} This assumption is fulfilled in many practical situations. 
For example, if we consider all modes with a decay rate  
less than unity and a frequency greater than $1$ Hz 
(i.e. $\alpha < 1$ and $\omega > 2\pi $), then this assumption is practically satisfied. \\

\noindent {\bf Assumption~3}. The eigenfrequencies $\omega_1$ and $\omega_2$ change much faster with $\gamma$ 
than the corresponding conventional PFs do, i.e., for each $k, l$:
 \begin{equation} \label{assump-3}
 \frac{1}{\omega_1}\frac{d\omega_1}{d\gamma}, \,  \frac{1}{\omega_2}\frac{d\omega_2}{d\gamma} 
 >> \frac{1}{p_{k1l}}\frac{d p_{k1l}}{d\gamma}, \, \frac{1}{p_{k2l}}\frac{d p_{k2l}}{d\gamma} \ .
\end{equation}

The characteristic property of LMIEs and pair LPFs is established by the following proposition.\\

\noindent {\bf Property~3}.
{\it Let the matrix ${\bf A}$ change under assumptions} 1, 2, 3 {\it in such a way that 
it passes through the point $\omega_1 = \omega_2$.
Then, the LMIE $\overline{E}_{z \,i j}$ in \eqref{def-3} and 
Lyapunov energies $E_{{x_k}12}(\gamma)$ in \eqref{Pxkij} associated with a pair of eigenvalues 
$\lambda_1$ and $\lambda_2$ 
reach the local maximum  
in the neighborhood of point $\omega_1 \approx \omega_2$,
both under conventional \eqref{ini-cond-convent} and spherically symmetrical \eqref{assume-spheric-symmetry}
initial conditions.}\\

\noindent {\bf Proof.} 
Choose a state variable $k$ such that $|p_{ki}| \ne 0$, $|p_{kj}| \ne 0$. 
Under conventional initial conditions \eqref{ini-cond-convent}, 
according to (\ref{prop-1}), we have
\begin{equation} \label{proof-prop3-1}
E_{{x_k}12}=E_{x_k} \tilde{e}_{k(12)} = \frac{1}{4} \sum_{i=1,\, 1^*} \sum_{j=2, \, 2^*} {\rm Re}\, \left\{ - \frac{p^*_{ki} p_{kj}}{\lambda^*_i+\lambda_j} \right\} \ .
\end{equation}
Let us introduce the notations $\alpha = \alpha_1+\alpha_2$, $\omega=\omega_1+\omega_2$, $\Delta\omega = \omega_2-\omega_1$. 
According to definition (\ref{PF-def}), we have 
$p_{k1^*}= p^*_{k1}$ and $p_{k2^*}= p^*_{k2}\,$.
Using this and combining the terms in (\ref{proof-prop3-1}), we obtain
\begin{equation} \label{proof-prop3-2}
E_{{x_k}12} =  \frac{1}{2}\frac{\alpha C^k_R+\Delta\omega C^k_I}{\alpha^2+\Delta\omega^2} 
+ \frac{1}{2}\frac{\alpha D^k_R+\omega D^k_I}{\alpha^2+\omega^2} \ , 
\end{equation}
where the coefficients $ C^k_R = {\rm Re}\,\{p_{k1} p^*_{k2} \}$, $C^k_I = {\rm Im}\,\{p_{k1} p^*_{k2}\}$, 
$D^k_R = {\rm Re}\,\{p_{k1} p_{k2} \}$, $D^k_I = {\rm Im}\,\{p_{k1} p_{k2}\}$ have the same order of magnitude.
According to assumption (\ref{assump-2}), $\alpha<< \omega$ 
and the second term in (\ref{proof-prop3-2}) is negligible  
\begin{gather} \label{proof-prop3-3}
\alpha \cdot E_{{x_k}12} =  \frac{1}{2}\frac{ C^k_R+\Delta\tilde{\omega}\, C^k_I}{1+\Delta\tilde{\omega}^2} + O(\epsilon) \ , \\
\text{where} \ \ \ \Delta\tilde{\omega} = \Delta\omega/\alpha, \ \ \epsilon=\alpha/\omega << 1 \ . \nonumber
\end{gather}
According to assumption (\ref{assump-3}), the dependence of the coefficients 
$C^k_R$ and $C^k_I$ on $\gamma$ can be neglected.  
Then, the maximum absolute value of \eqref{proof-prop3-3} is reached approximately at
$ C_I \Delta \tilde{\omega} \approx \sqrt {C_R^2+C_I^2} - C_R $, or, returning to the original notation, at 
\begin{gather}
|\Delta\omega| \approx \alpha \cfrac {|p_{k1} p^*_{k2}| - Re\{p_{k1} p^*_{k2}\}}{|{\rm Im}\,\{p_{k1} p^*_{k2}\}|} < \alpha_1+\alpha_2<< \omega_1+\omega_2 \ . \nonumber
\end{gather}
Under spherically symmetrical initial conditions \eqref{assume-spheric-symmetry}, 
the proof is similar, albeit with the replacement of $C^k_R$ and $C^k_I$ by 
\[
\widetilde{C}^k_R= \sum^n_{l=1} {\rm Re}\,\{p_{k1l} p^*_{k2l}\} \ \ \ \text{and} \ \ \ \widetilde{C}^k_I = \sum^n_{l=1} {\rm Im}\,\{p_{k1} p^*_{k2}\}. 
\]
According to (\ref{Ezij=sumExkij}), $\overline{E}_{z12} = \sum^n_{k=1} E_{{x_k}12} $.
Therefore, under spherically symmetrical initial conditions \eqref{assume-spheric-symmetry}, the function $\overline{E}_{z12}$
has the same structure as (\ref{proof-prop3-3}):  
\begin{equation} \label{proof-prop3-4}
\overline{E}_{z12} = \sum^n_{k=1} E_{{x_k}12} 
\approx 
\frac{1}{2\alpha} \frac{(\sum_{k} \widetilde{C}^{k}_R)+\Delta\tilde{\omega} \, (\sum_{k} \widetilde{C}^{k}_I)}{1+\Delta\tilde{\omega}^2} \ . \nonumber
\end{equation}
This function reaches its local maximum approximately at
\begin{equation}
\Delta\omega \approx \alpha \left(\sqrt{\sigma^2+1}-\sigma\right)< \alpha_1+\alpha_2 << \omega_1+\omega_2 \ , \nonumber
\end{equation}
where $\sigma = \sum_{k}  \widetilde{C}^{k}_R \Big/ \sum_{k}  \widetilde{C}^{k}_I $. 
$\square$ \\

\noindent {\bf Remark.} 
It follows from Property~3 that if the Lyapunov energies corresponding to other mode pairs
 change rather slowly in the neighborhood of point $\omega_1 \approx \omega_2$,  
 then the functions $LMIF_{12}(\gamma)$ in \eqref{def-4} and MISLPF  $\tilde{e}_{k(12)}(\gamma)$ in \eqref{def-5-1}
 also reach a local maximum in the neighborhood of point $\omega_1 \approx \omega_2$.\\

It follows from Property~3 that the values of pair Lyapunov energies and pair LPFs facilitate the identification of the resonant
interactions between lightly damped oscillating modes. In general, it can be seen that 
in \eqref{prop-1}, the value of LPF $\epsilon_{ki}=\sum_j \tilde{\epsilon}_{kij}$ is not simply
proportional to $|p_{ki}|^2$, but involves the interaction of the $i$-mode
with other $j$ modes. In accordance with Property~3, 
the closer the $j$-mode is in frequency to the $i$-mode, 
the greater its contribution $\tilde{\epsilon}_{kij}$ to $\epsilon_{ki}$ is, 
which characterizes {\it the resonance interaction}.
Similarly,  the smaller $|Re(\lambda_j)|$ is, 
the greater the $j$-mode contribution $\tilde{\epsilon}_{kij}$ to $\epsilon_{ki}$. 
This characterizes {\it the interaction with a weakly stable mode}. 

\subsection{Fast calculation of Lyapunov PFs}

The possibility of fast calculation of Lyapunov energies and LPFs for the purposes of selective modal analysis 
is based on formulas \eqref{Pxki}, \eqref{MISLPF-ini-cond-1}, \eqref{eki-evx}, \eqref{Pxkij}, \eqref{pair-MISLPF-ini-cond-1}, \eqref{ekij-sym-ini}. 
These can be summarized as follows. \\

\noindent {\bf Property~4}.
{\it The Lyapunov energies of states associated with particular modes and unnormalized LPFs 
can be calculated through the corresponding eigenvalues and eigenvectors, namely 
under conventional initial conditions \eqref{ini-cond-convent}}:
\begin{gather}
E_{x_k i} = E_{x_k} e_{ki} = - \left\{ (u^k_i)^* (\varv^k_i)^* \ {\bf e}^T_k (\lambda^*_i {\bf I} +{\bf A})^{-1} {\bf e}_k \right\}_H \ , \nonumber \\
E_{x_k i j} = E_{x_k} \tilde{e}_{kij} = -  \frac{1}{2} \left\{ \cfrac{(\varv^k_i)^* (u^k_i)^* u^k_j \, {\varv}^k_j}{\lambda_i^* +\lambda_j}
+ \cfrac{{\varv}^k_i \, u^k_i \, u^k_j \, {\varv}^k_j}{\lambda_i +\lambda_j} \right\}_H \ , \nonumber
\end{gather}
{\it and under spherically symmetric initial conditions \eqref{assume-spheric-symmetry}}:
\begin{gather}
E_{x_k i} = E_{x_k} e_{ki} = - {\rm Re} \left\{ \mathrm{trace} \left( ({\boldsymbol{\varv}}^T_i)^* (u^k_i)^*  
{\bf e}^T_k (\lambda^*_i {\bf I} +{\bf A})^{-1} \right) \right\} \ , \nonumber \\
E_{x_k i j} = E_{x_k} \tilde{e}_{kij} = -  \frac{1}{2} \mathrm{trace}
\left\{ \cfrac{(\boldsymbol{\varv}^T_i)^* (u^k_i)^* u^k_j \, \boldsymbol{\varv}^T_j}{\lambda_i^* +\lambda_j}
+ \cfrac{\boldsymbol{\varv}_i \, u^k_i \, u^k_j \, \boldsymbol{\varv}^T_j}{\lambda_i +\lambda_j} \right\}_H \ . \nonumber
\end{gather}
Property~4 implies that in order to analyze the dynamics of Lyapunov energies and unnormalized LPFs 
in the critical part of the spectrum, it is not necessary to know the entire spectrum of the system matrix. 
It is sufficient to know only the eigenvalues and eigenvectors in the critical part of the spectrum, which can be obtained 
using the well-known algorithms of selective modal analysis, 
such as the modified Arnoldi method or simultaneous iterations \citep{wang:90}. 
Therefore, the LMA can be performed quickly for the critical modes and 
serve as a basis for the fast calculation of their behavior in large-scale dynamical systems in real time.\\

\noindent {\bf Remark}. 
Formulas of Proposition~4 are valid only for systems with a simple spectrum. 
When two eigenvalues approach each other, these formulas start to become ill-conditioned.
When a multiple root appears, the eigenvectors may not be uniquely determined. 
In this case, all the formulas must be modified. However, this is beyond the scope of this work. 

\subsection{LPFs and eigenvalue sensitivities}

\noindent {\bf Property~5}.
 {\it LPFs and corresponding Lyapunov energies can be represented 
using the sensitivities of the corresponding eigenvalues with respect to the elements of 
the dynamic matrix ${\bf A}$}: 
\begin{gather}  
E_{x_k} \, e_{ki} = - \left\{ \left( \frac{\partial\lambda_i}{\partial a_{kk}}  \right)^*{\bf e}^T_k (\lambda^*_i {\bf I} +{\bf A})^{-1} {\bf e}_k \right\}_H \ , \nonumber \\
E_{x_k} \, \tilde{e}_{k(ij)} = - \frac{1}{2}
\left\{ \frac{1}{\lambda^*_i+\lambda_j} \left( \frac{\partial\lambda_i}{\partial a_{kk}}  \right)^* \frac{\partial\lambda_j}{\partial a_{kk}}
+ \frac{1}{\lambda_i+\lambda_j} \frac{\partial\lambda_i}{\partial a_{kk}}  
\frac{\partial\lambda_j}{\partial a_{kk}}  
 \right\}_H . \nonumber  
\end{gather}

\noindent These formulas are obtained by substitution of (\ref{PF-sensetivities}) into 
the formulas of Proposition~1. 
In practical applications, the corresponding sensitivity coefficients can be identified using various measurement methods.\\

Based on the above considerations,
it can be expected that the LPFs may provide additional
advantages for the purposes of the small-signal stability
analysis, while the PFs retain their importance in assessing
the instantaneous dynamics of the system. \\

\section{Simulation experiment} 

In this section, we present a fairly simple numerical experiment,
which is nevertheless able to demonstrate the potential advantages of Lyapunov modal analysis.

\subsection{Experiment description}

\begin{figure}
  \centerline{
  \includegraphics[width=8.9cm]{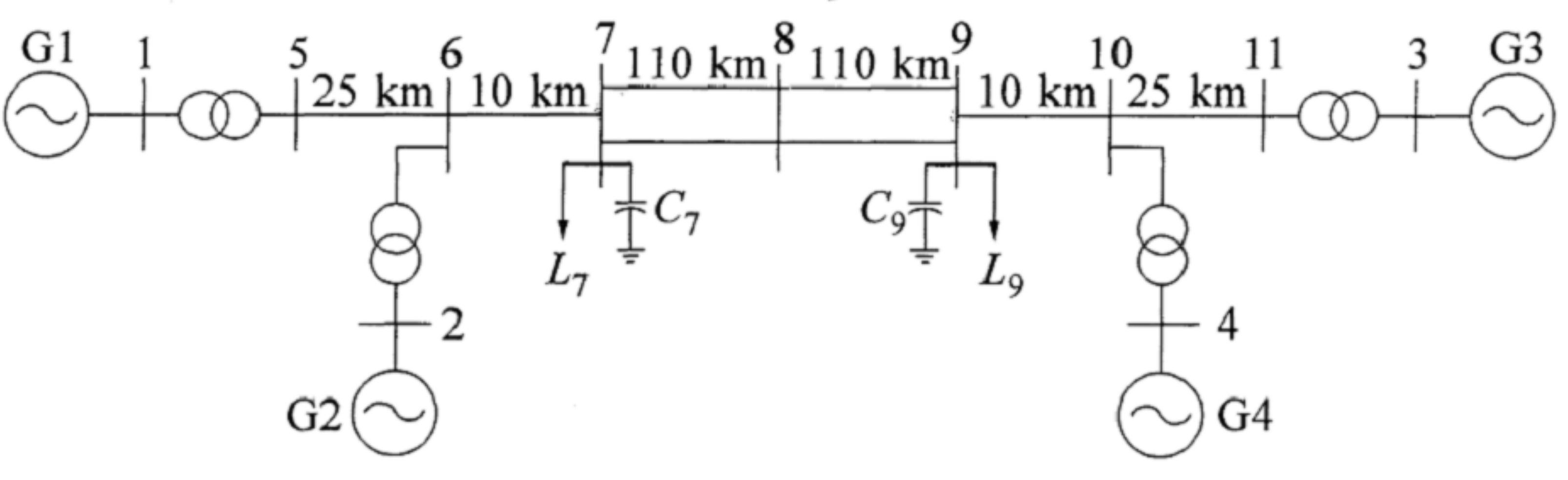}
   }
 \caption{Two-area power system with four generators from (\cite{kundur:94}).}
\label{fig:kundur}
\end{figure}

For carrying out a simulation experiment, we use the two-area power system with four generators that was considered by \citet{kundur:94} and is shown in Figure~\ref{fig:kundur}. It consists of two similar areas connected by a weak tie. Each area consists of two coupled stations. The third station, G3, was considered as a swing bus. For dynamical analysis, all generators of the power system are represented by sixth order models. The speed controllers were ignored. The model parameters chosen were the same as those in the textbook (\citet{kundur:94}, Example 12.6, p.813). To avoid zero eigenvalues in the dynamic matrix, all rotor angles and speed deviations were taken with respect to those of reference generator G3.

\begin{table}
	\centering
	\begin{tabular}{ ccp{4.5cm}} 
		\hline
		Mode & Initial eigenvalue & Type and location \\ 
		\hline
$S_1$ & $-0.096$ & aperiodic rotor angle mode (mostly between G3 and G4) \\ \hline
$S_2$ & $-0.117$ & aperiodic rotor angle mode (mostly between G1 and G2) \\ \hline
$S_3$ & $-0.111\pm 3.43j$ & inter-area rotor oscillation (between G1, G2, G3, G4) \\ \hline
$S_4$ & $-0.265$ &local inter-machine flux linkage mode (between G3 and G4) \\ \hline
$S_5$ & $-0.276$ & local inter-machine flux linkage mode (between G1 and G2)  \\ \hline
$S_6$ & $-0.492\pm 6.82j$ & local inter-machine oscillation (between G1 and G2) \\ \hline
$S_7$ & $-0.506\pm 7.02j$ & local inter-machine oscillation (between G3 and G4) \\ \hline
		\hline
	\end{tabular}
	\caption{Modes and initial eigenvalues.}
	\label{tab:modes}
\end{table}

\begin{figure*} [t]
  \centerline{
  \includegraphics[width = 0.45\textwidth]{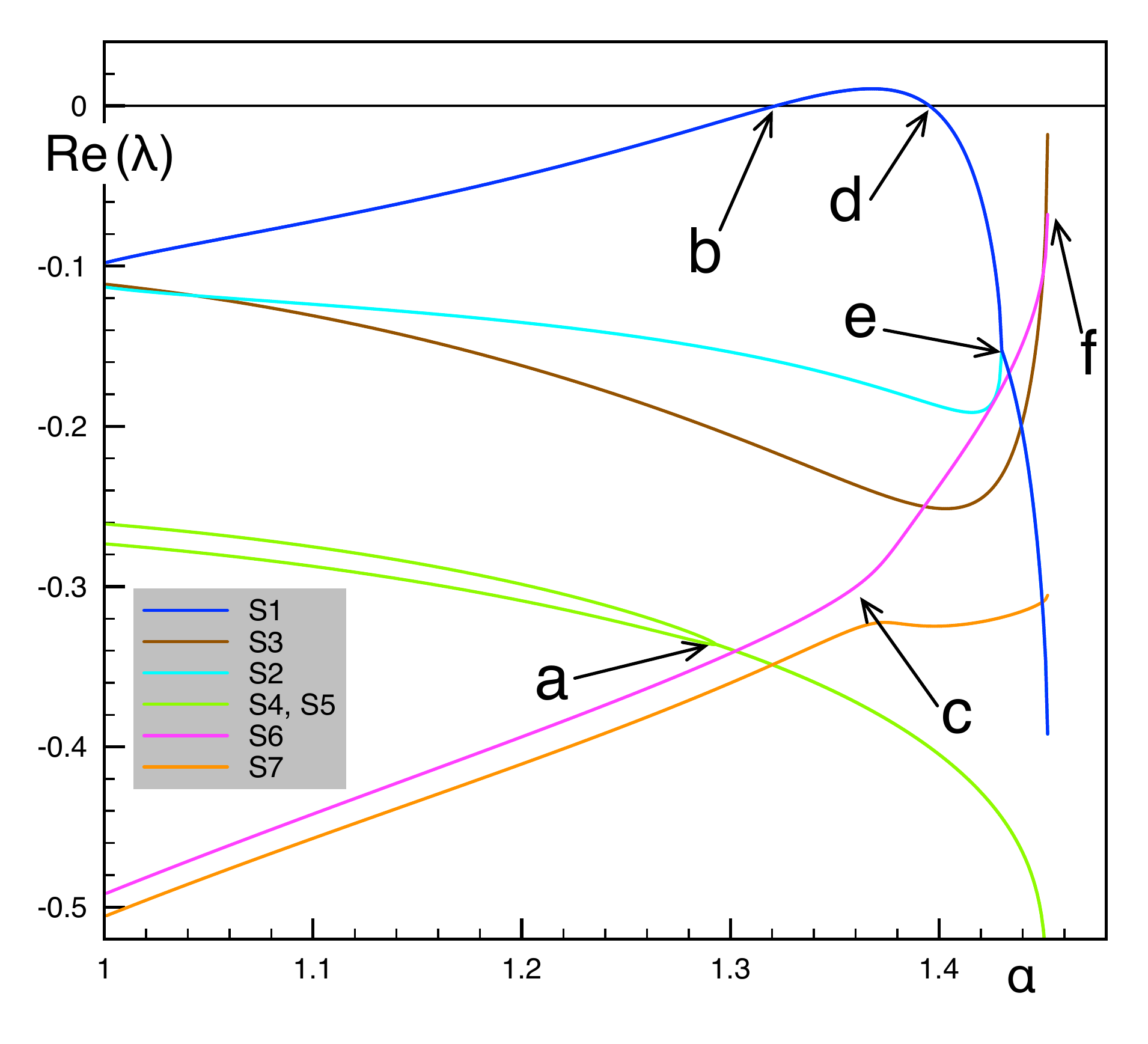} \ \ \ \ \ \  
  \includegraphics[width = 0.45\textwidth]{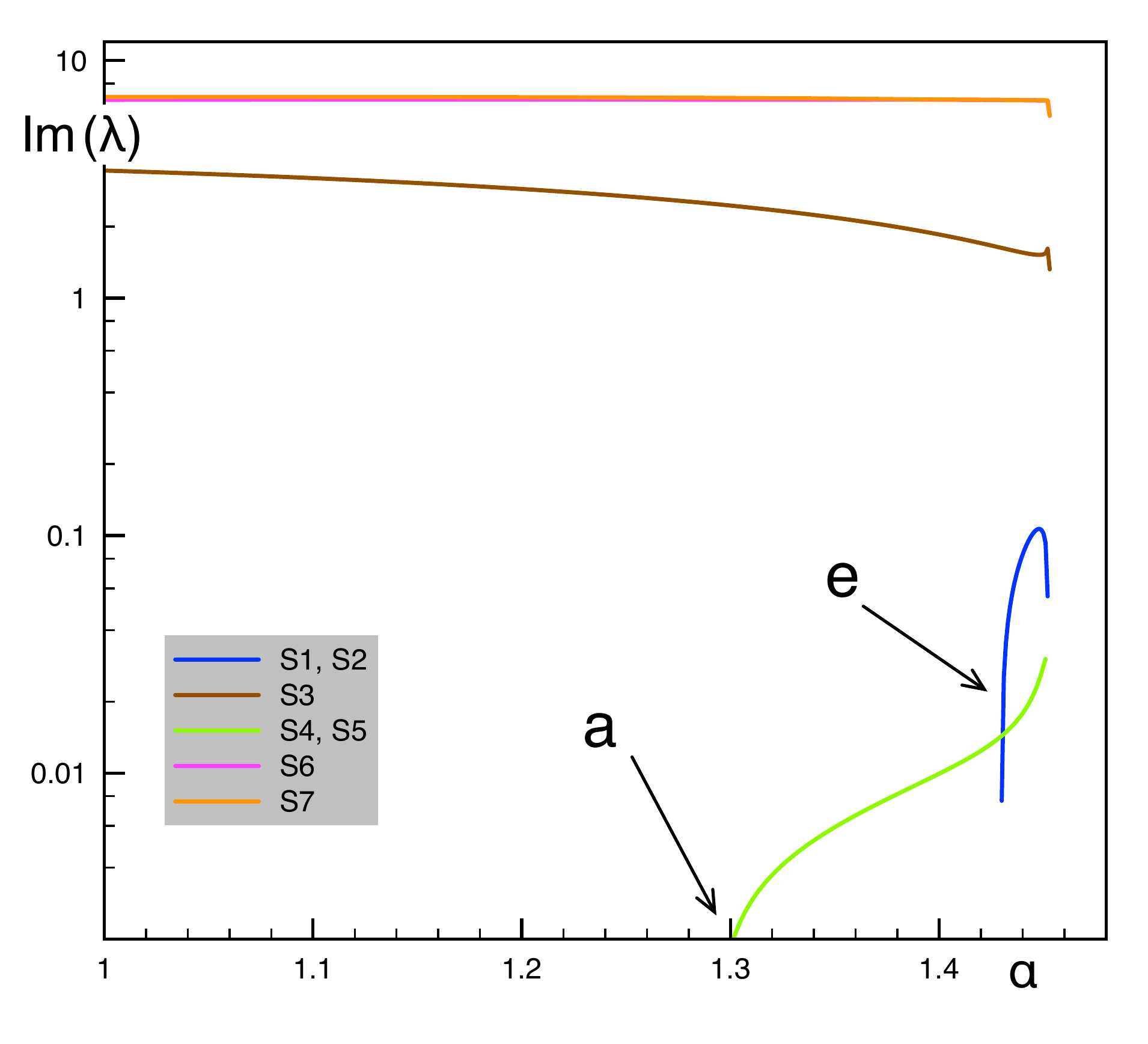} 
  }
   \caption{
   Trajectory of the real (on the left) and imaginary (on the right) parts of eigenvalues in the experiment. 
    }
   \label{fig:eigenvalues}
\end{figure*}

\begin{figure*} [!b]
  \centerline{
  \includegraphics[width = 0.45\textwidth]{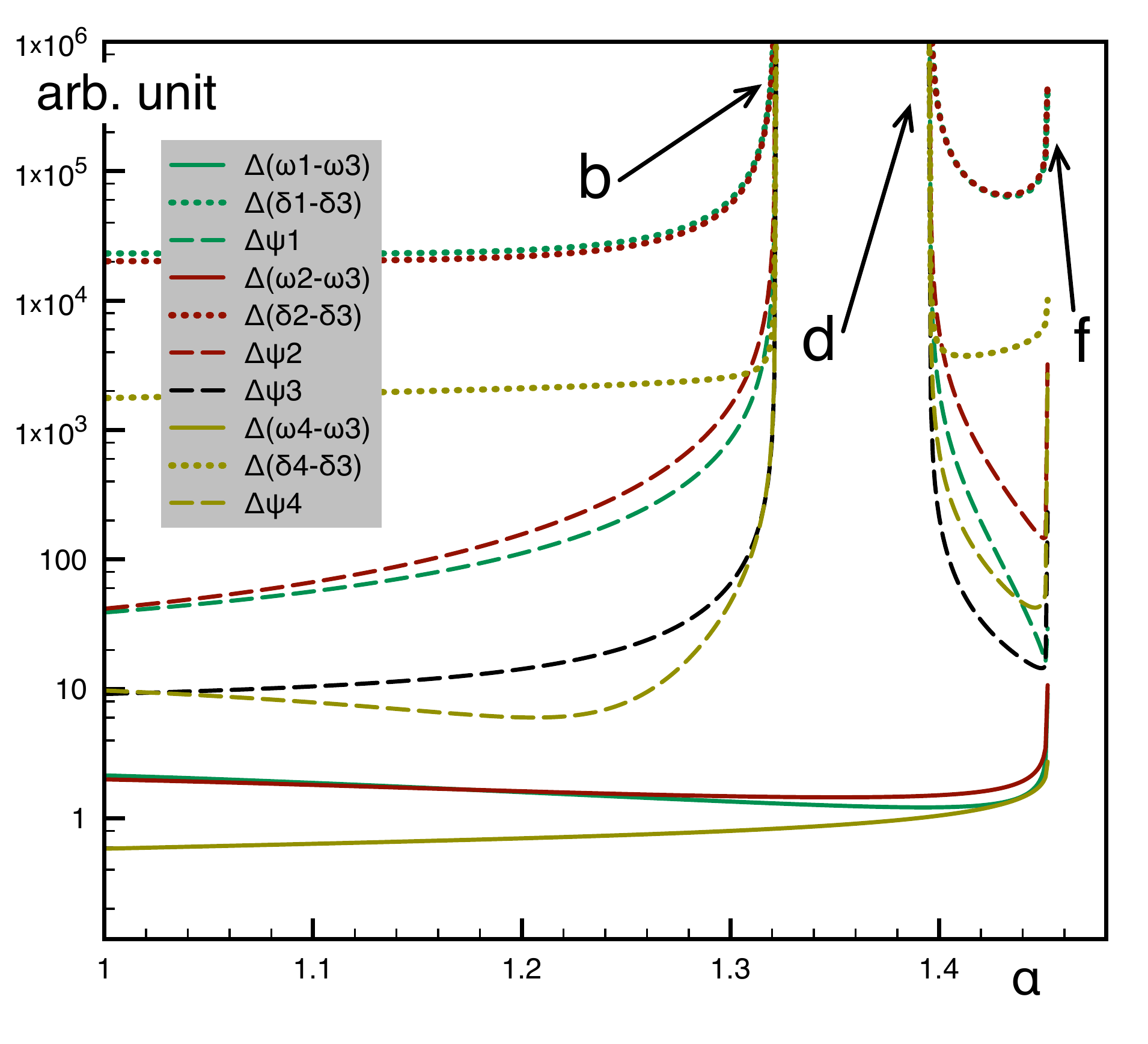} \ \ \ \ \ \  
  \includegraphics[width = 0.45\textwidth]{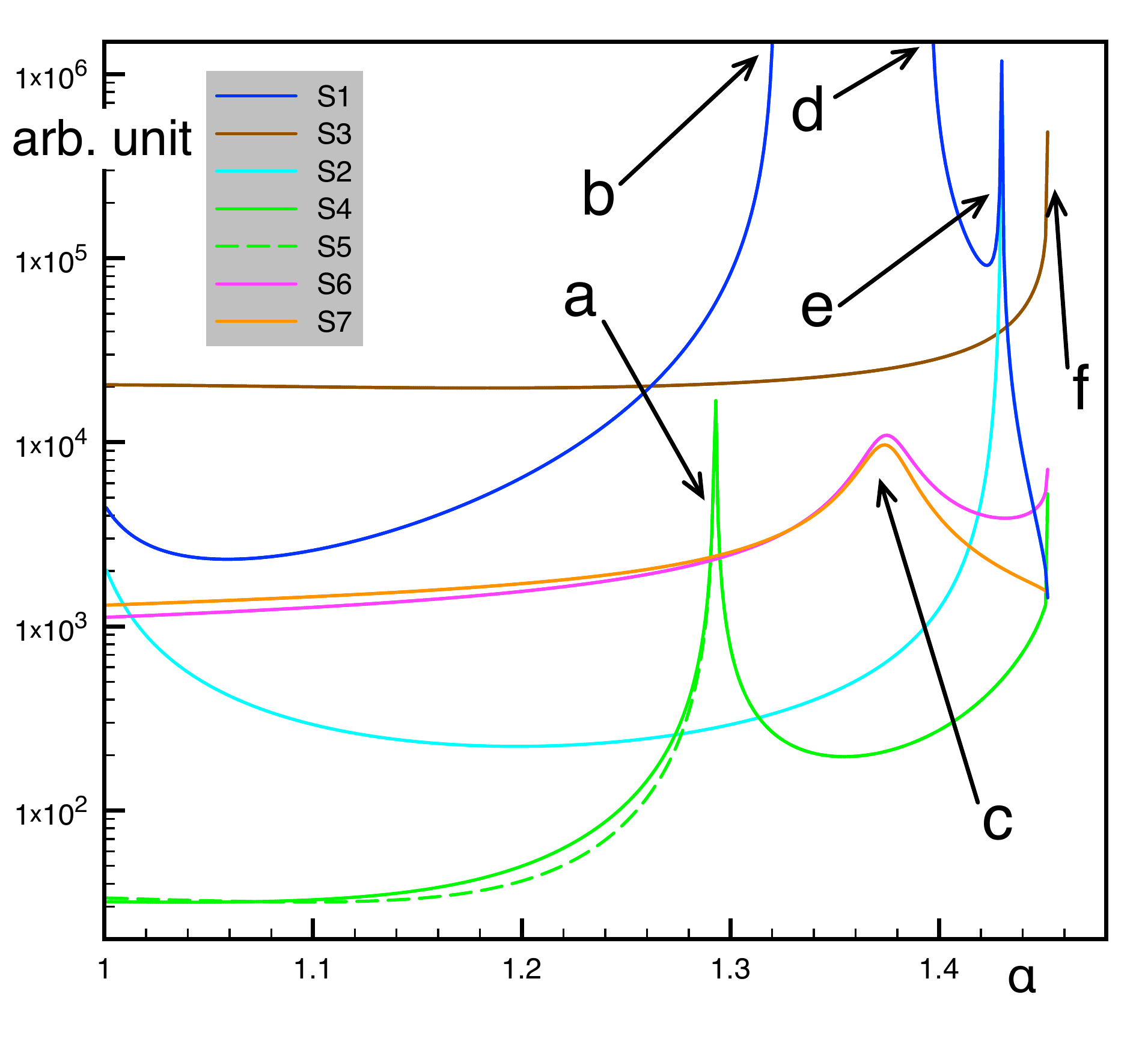} 
  }
   \caption{
   Evolution of Lyapunov energies of the states (on the left) and modes (on the right) in the experiment. 
    }
   \label{fig:lyap-energies}
\end{figure*}

We studied the limit of the system stability by simultaneously increasing all the loads ($P_L,Q_L,Q_C$) and the active power of each generator ($P_i$) while keeping the ratios between them fixed. We define the power increase coefficient as $\alpha =P/P_0=Q/Q_0$.
This change in $\alpha$ results in further changes in system modes. The list of the main ones is presented in Table~\ref{tab:modes}. All oscillations listed in the table are rotor angle electro-mechanical oscillations. The aperiodic modes S1 and S2 correspond to the rotor angle. The aperiodic modes S4 and S5 correspond to flux linkage.
Figure~\ref{fig:eigenvalues} shows the evolution of the real (on the left) and imaginary (on the right) parts of the eigenvalues during the test experiment 
as functions of the power increase coefficient $\alpha$. Mode names on the legend correspond to Table~\ref{tab:modes}.
The simulation indicates nontrivial dynamics of system modes. 
As the parameter $\alpha$ increases, the following changes, which are marked in Figure~\ref{fig:eigenvalues}, are observed:

\begin{enumerate}[(a)]
\item At $ \alpha_{\rm a} \approx 1.293$, the aperiodic S4 and S5 modes merge into one low-frequency oscillation.
\item When $ \alpha_{\rm b} \approx 1.321$, the aperiodic inter-area angle mode S1 becomes unstable. As can be seen, however, the instability of the S1 mode does not influence the behavior of the other modes.
\item When $ \alpha_{\rm c} \approx 1.375$, there is resonance between the S6 and S7 oscillations, after which they become inter-area oscillations.
\item The S1 mode becomes stable again at $\alpha_{\rm d} \approx 1.395$.
\item When $ \alpha_{\rm e} \approx 1.430$, the aperiodic modes S1 and S2 merge into one low-frequency oscillation.
\item The system becomes unstable again at $\alpha_{\rm f} \approx 1.453$. An obvious relationship 
can be seen in the behavior of the dangerous modes S3, S6, and S7 in the pre-fault operation.
\end{enumerate}

\begin{figure*} [!ht]
  \centerline{
  \includegraphics[width = 0.4\textwidth]{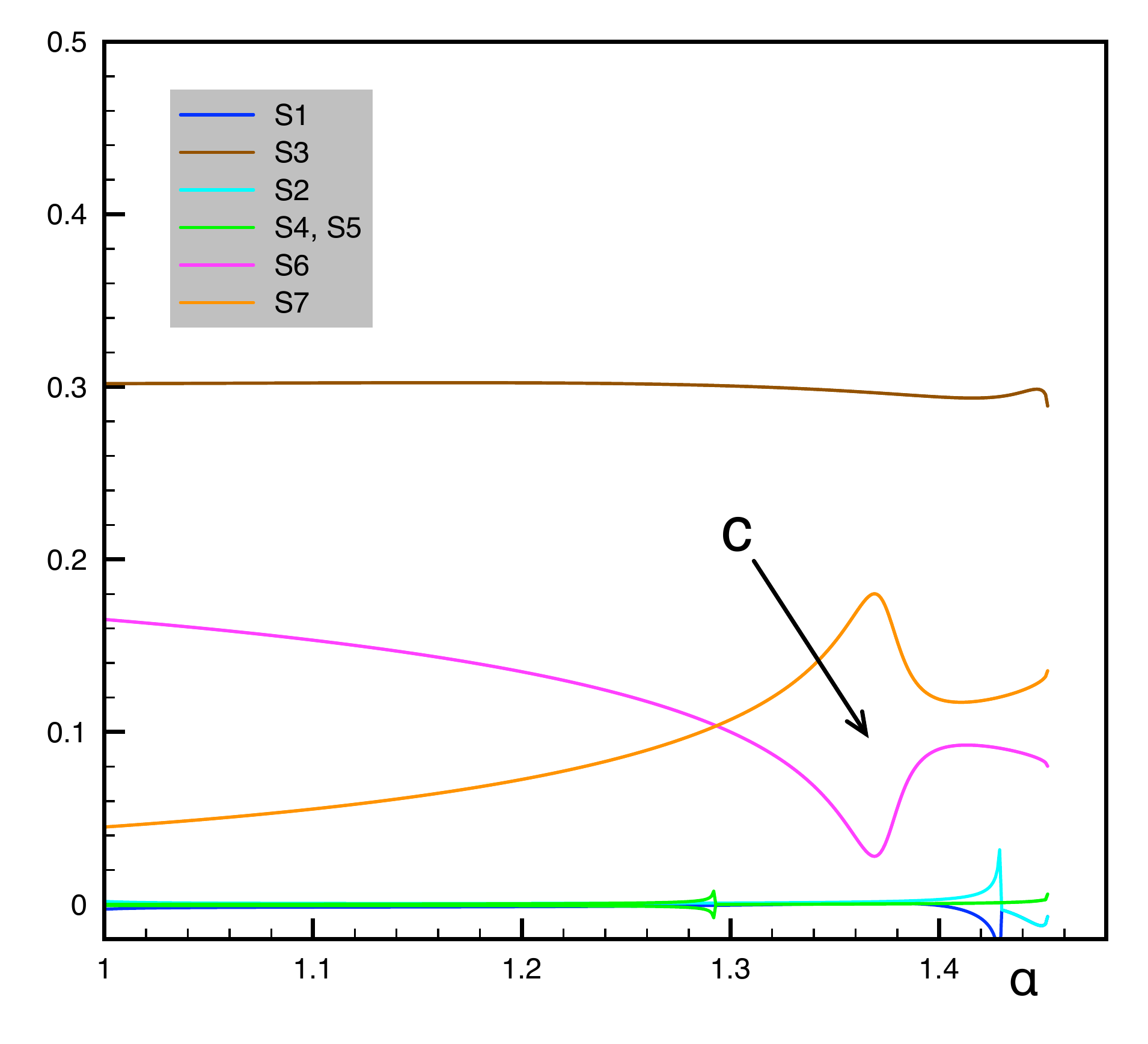} \ \ \ \ \ \  
  \includegraphics[width = 0.4\textwidth]{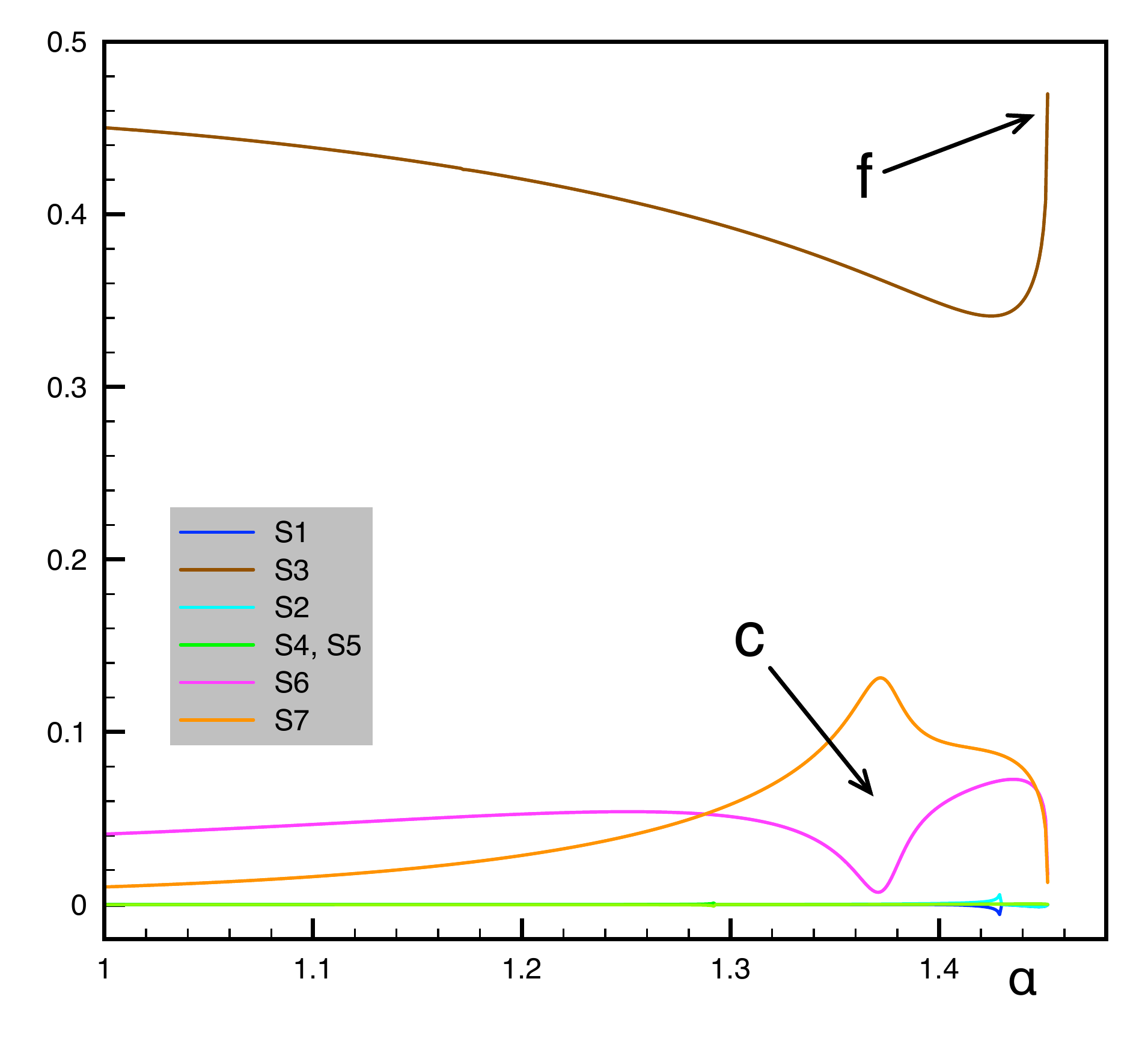} 
  }
  \centerline{
  \includegraphics[width = 0.4\textwidth]{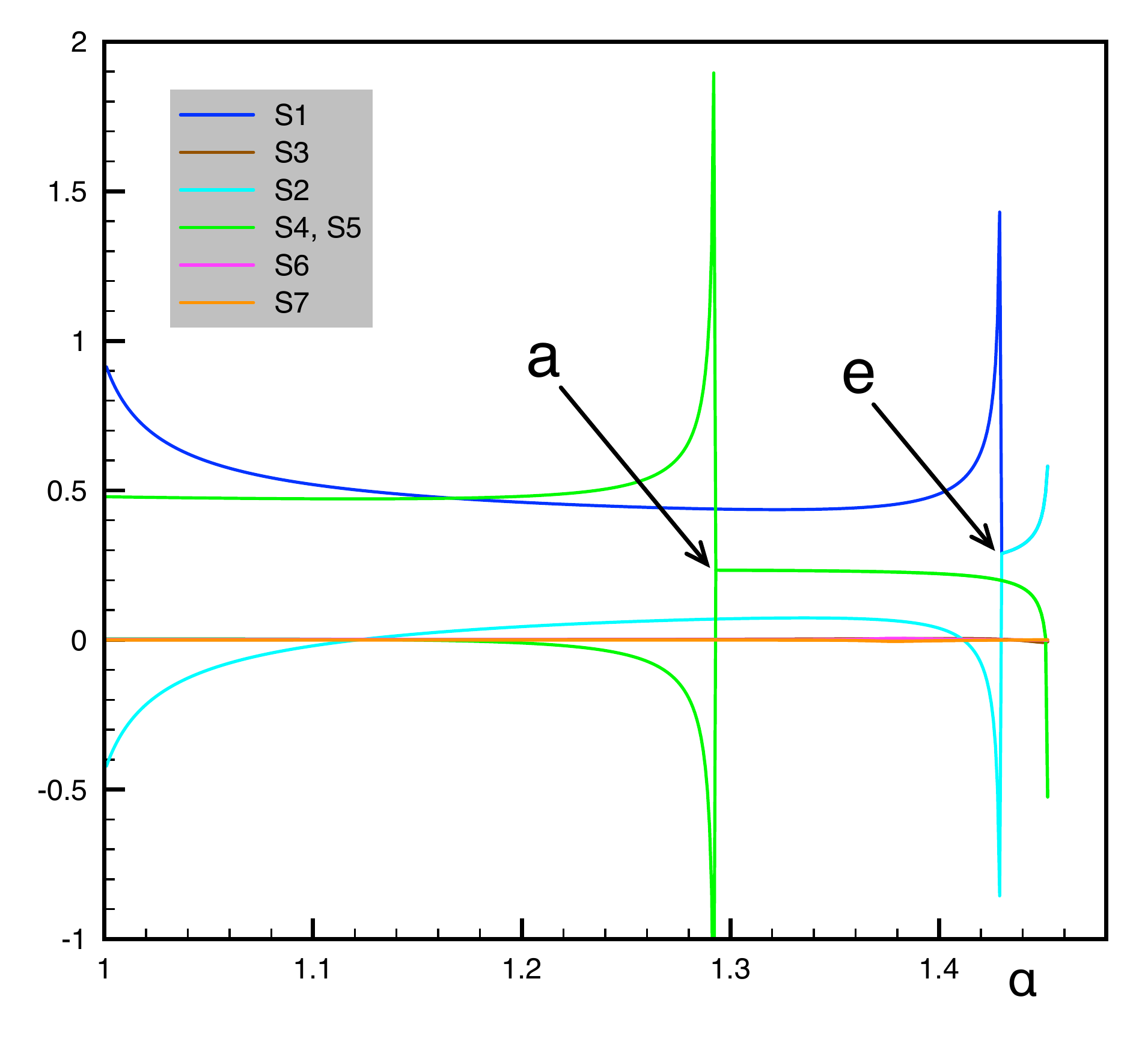} \ \ \ \ \ \  
  \includegraphics[width = 0.4\textwidth]{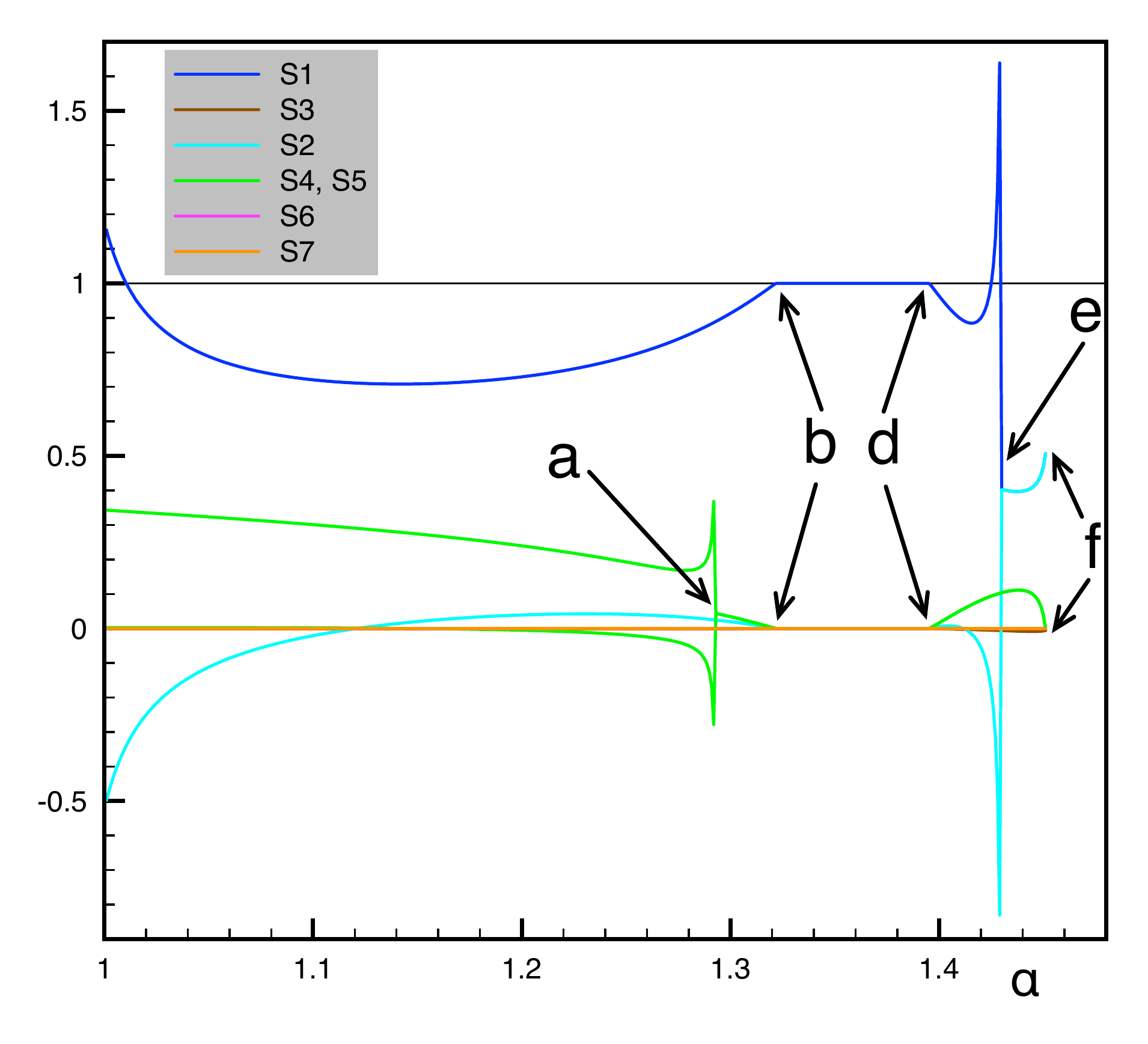} 
  }
   \caption{
Behavior of the conventional PFs $p_{ki}$ (on the left) and Lyapunov PFs $e_{ki}$ (on the right) for the state variables
$\Delta(\omega_1-\omega_3)$ and $\Delta\Psi_1$ above and below, respectively.
    }
   \label{fig:MISPFs}
\end{figure*}

\subsection{Simulation results and discussion}

The evolution of Lyapunov energies of the states and modes is shown in Figure~\ref{fig:lyap-energies}
as a function of the power increase coefficient $\alpha$.
The units of Lyapunov energies are derived from the system of units used in (\cite{kundur:94}, see also Remark~2). 
The graphs on the left show the values of $E_{x_k}$ calculated by \eqref{Lener-k-state} for the main state variables, namely for 
the deviations of flux linkages and rotor angles and speeds of different generators. 
When the system loses stability at points (b), (d) and (f), the Lyapunov energies in the corresponding variables tend to infinity.
However, in some other variables, they remain bounded.
The graphs on the right show the values of $|\boldsymbol{u}_i|^2 E_{z_i}$, which are calculated using \eqref{Lener-i-mode+} and 
characterize the invariant measure of Lyapunov modal energies (see Remark~1) for the modes $S_1$ to $S_7$ listed in Table~\ref{tab:modes}.
Each curve characterizing an oscillation contains two identical graphs; either of which corresponds to one of
the complex conjugate eigenvalues.
The graphs reflect all qualitative changes in the spectrum of the system, including the loss of stability by the corresponding modes at points (b), (d), and (f), the fusion of aperiodic modes at points (a) and (e), and the resonance between the oscillations $S_6$ and $S_7$ at (c).

Figure~\ref{fig:MISPFs} shows the behavior of the conventional MISPF $p_{ki}$ defined by \eqref{PF-def} (on the left) and MISLPF $e_{ki}$ defined by \eqref{def-1} and \eqref{MISLPF-ini-cond-1} (on the right) depending on $\alpha$. The plots at the top show the modal PFs in $\Delta(\omega_1-\omega_3)$, i.e., in the deviation of the rotor angle speed of generator G1 with respect to that of generator G3. The plots at the bottom show the modal PFs in $\Delta\Psi_1$, i.e., in the deviation of flux-linkage of generator G1. The general composition of modes in the considered state variables according to both $p_{ki}$ and $e_{ki}$ are similar.
Both coefficients identify the process of merging the aperiodic modes at points (a) and (e).
However, unlike conventional PFs, Lyapunov PFs clearly identify the moments of stability loss occurring due to the corresponding mode and state variables. 
In the plots at the top, two graphs of $e_{ki}$ characterizing the oscillation $S_3$, which loses stability at (f), in the sum tend to unity at $\alpha_f$. 
In the plots at the bottom, the graph characterizing the aperiodic mode $S_1$, which loses stability in the interval between points (b) and (d), tends to unity in the same interval. Thus, in accordance with Property~2, the MISLPFs identify the stability loss occurring due to specific modes and state variables.

\begin{figure*} 
  \centerline{
  \includegraphics[width = 0.4\textwidth]{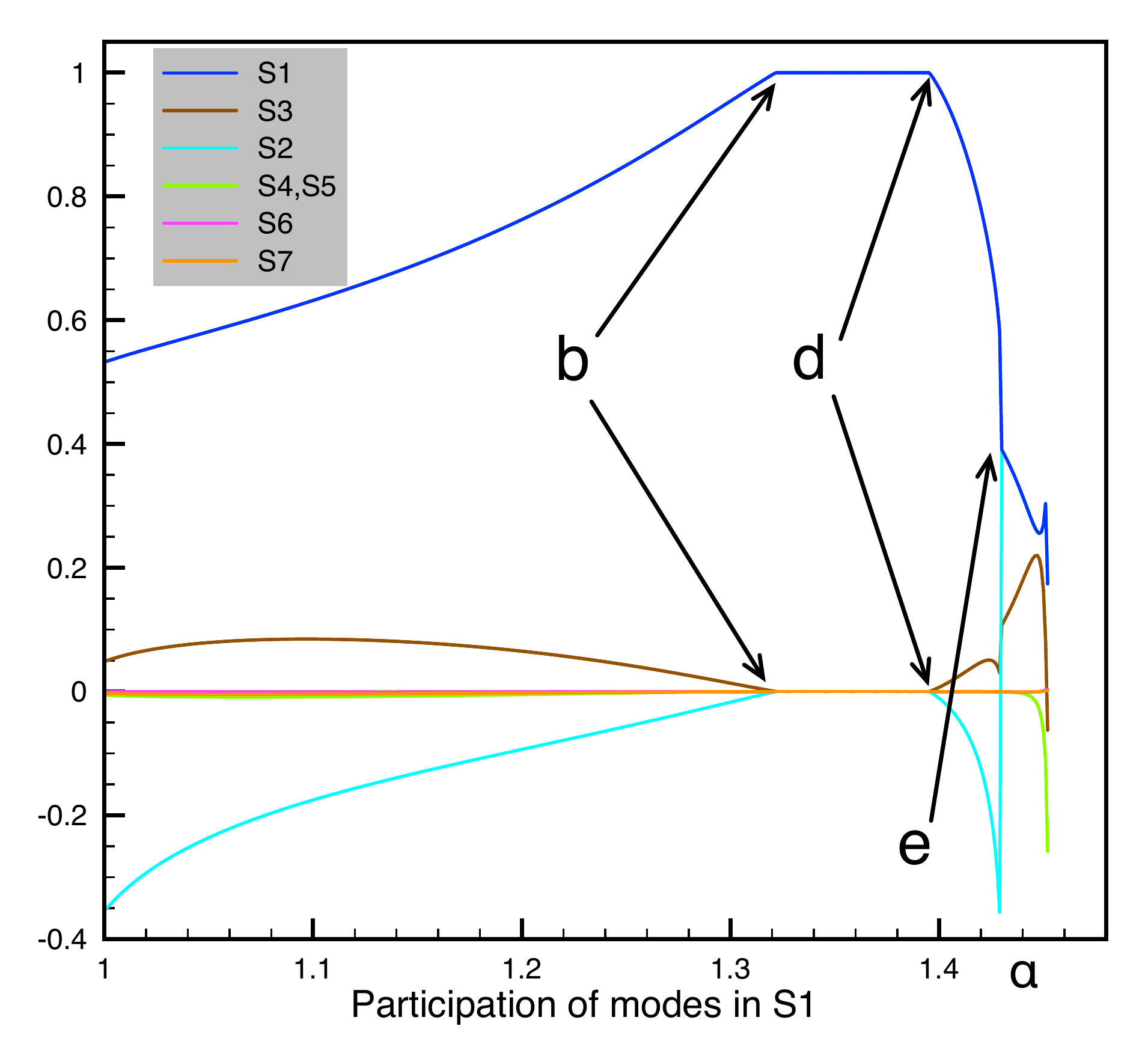} \ \ \ \ \ \  
  \includegraphics[width = 0.4\textwidth]{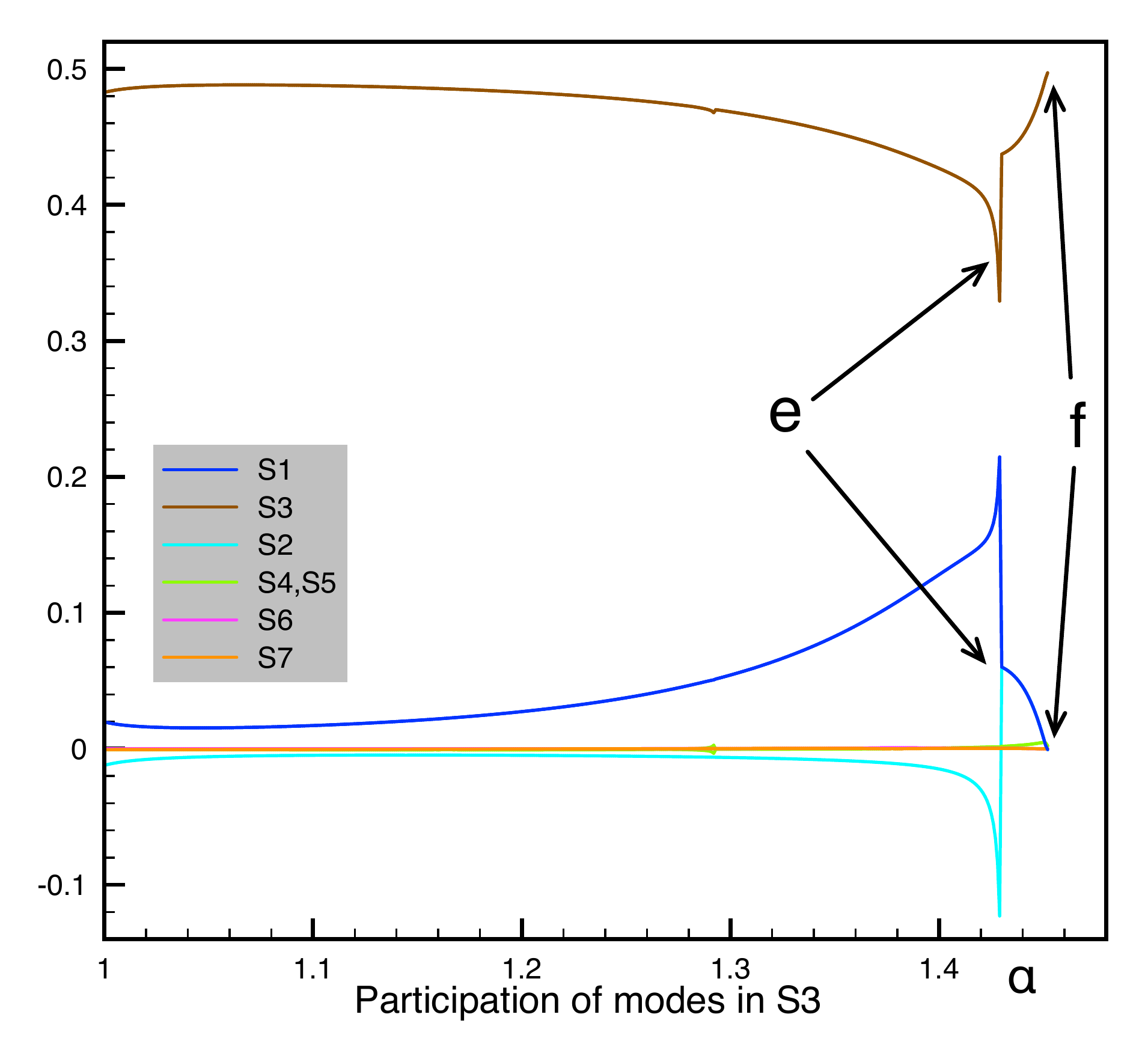} 
  }
  \centerline{
  \includegraphics[width = 0.4\textwidth]{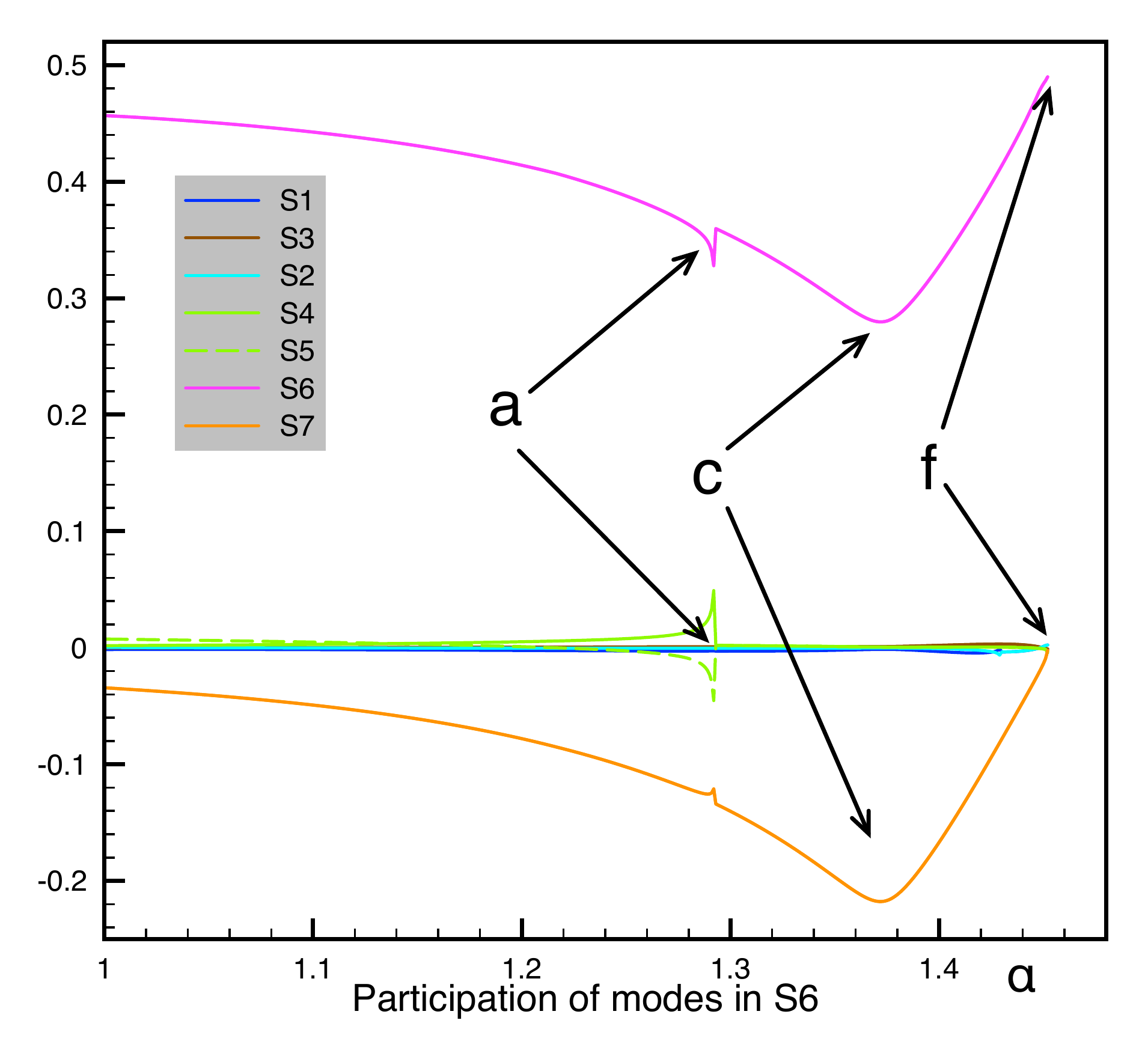} \ \ \ \ \ \  
  \includegraphics[width = 0.4\textwidth]{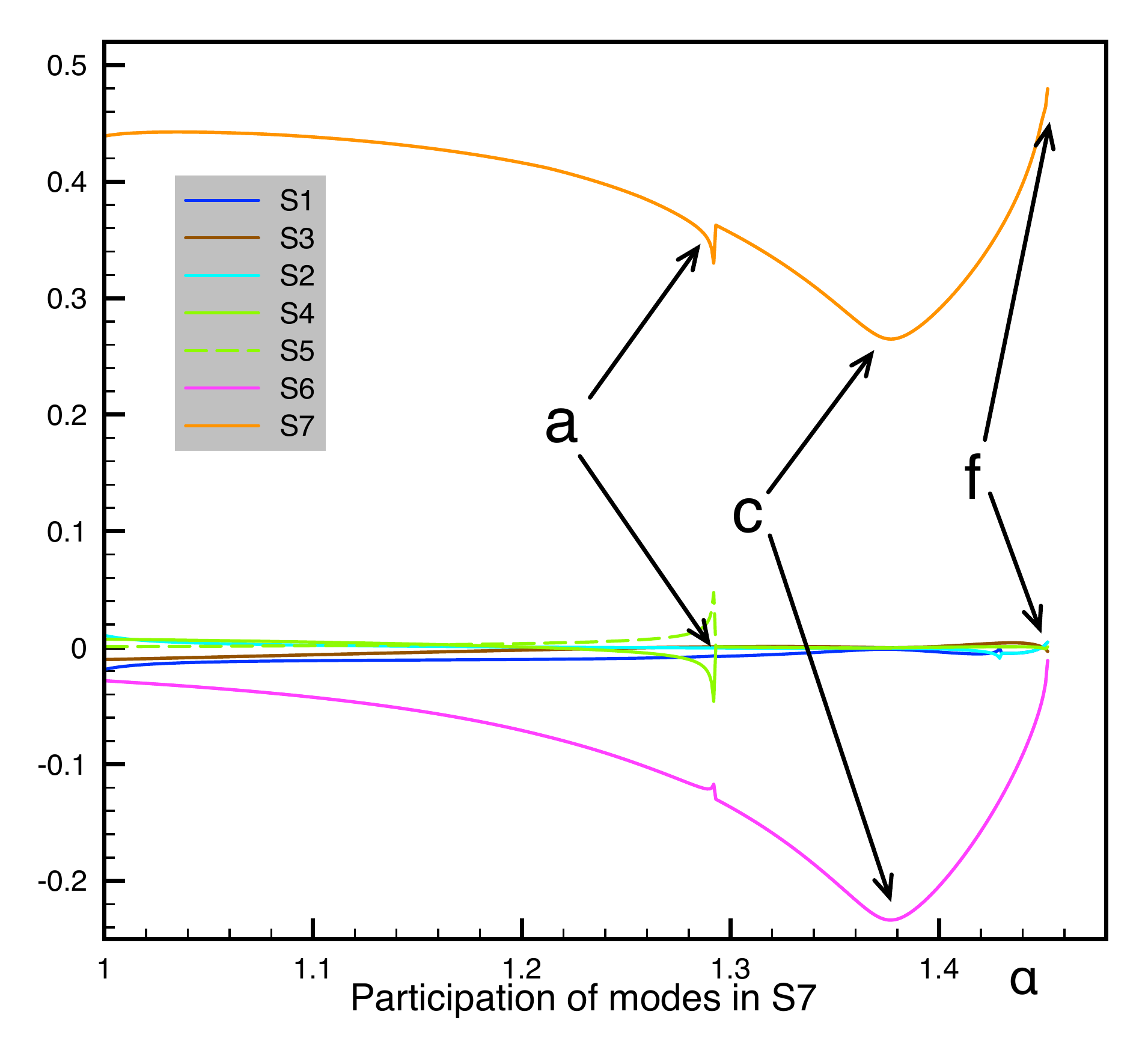} 
  }
   \caption{
Lyapunov modal interaction factors normalized to 1.
    }
   \label{fig:LMIFs}
\end{figure*}

Figure~\ref{fig:LMIFs} shows the behavior of LMIFs calculated by \eqref{def-4} depending on $\alpha$.
The participations of all modes in modes $S_1$, $S_3$, $S_6$, and $S_7$ are shown on separate tabs. 
In the case of unstable modes, in accordance with the physical meaning, we assume that
\begin{equation} \nonumber
{\rm LMIF}_{ij} = +\infty \ , \ \ \text{if} \ \ {\rm Re} \{\lambda^*_i +\lambda_j\} > 0 \ .
\end{equation}
Each curve
showing the interaction with the oscillation contains two identical graphs, either of which corresponds to one of
the complex conjugate eigenvalues. In accordance with the chosen normalization, the sum of the absolute values
of the graphs on each tab is always equal to one. 
LMIF plots allow for the observation of a general structure of the
modal interaction, as well as to identify the following
characteristic features of the modal dynamics.
\begin{itemize}
\item 
{\it Loss of stability of an aperiodic or oscillatory mode.}
When the $S_1$ mode becomes unstable at point (b), its own participation approaches 1, and the
participations of the other modes in it disappear.
Similarly, when the oscillatory modes $S_3$, $S_6$, and $S_7$ approach the stability boundary
at (f), their own participations also tend to unity.

\item
{\it Merging of two aperiodic modes into one oscillation.}
When aperiodic modes $S_4$ and $S_5$ merge into a
single oscillation at (a), there is a noticeable increase in the participations of merging
modes in other modes before merging and a sharp increase in the participations of other modes after that.
A similar phenomenon is observed when aperiodic modes $S_1$ and $S_2$ merge into a single
oscillation at (e). 

\item 
{\it Occurrence of a resonance between two oscillations.}
Oscillations $S_6$ and $S_7$ interact
mainly with each other (see the tabs at the bottom of Figure~\ref{fig:LMIFs}). As their frequencies
approach each other at (c), the graphs show a characteristic increase in the mutual participation
of these modes in each other, with the opposite sign.
\end{itemize}

\begin{figure} [ht]
  \centerline{
  \includegraphics[width=7.8cm]{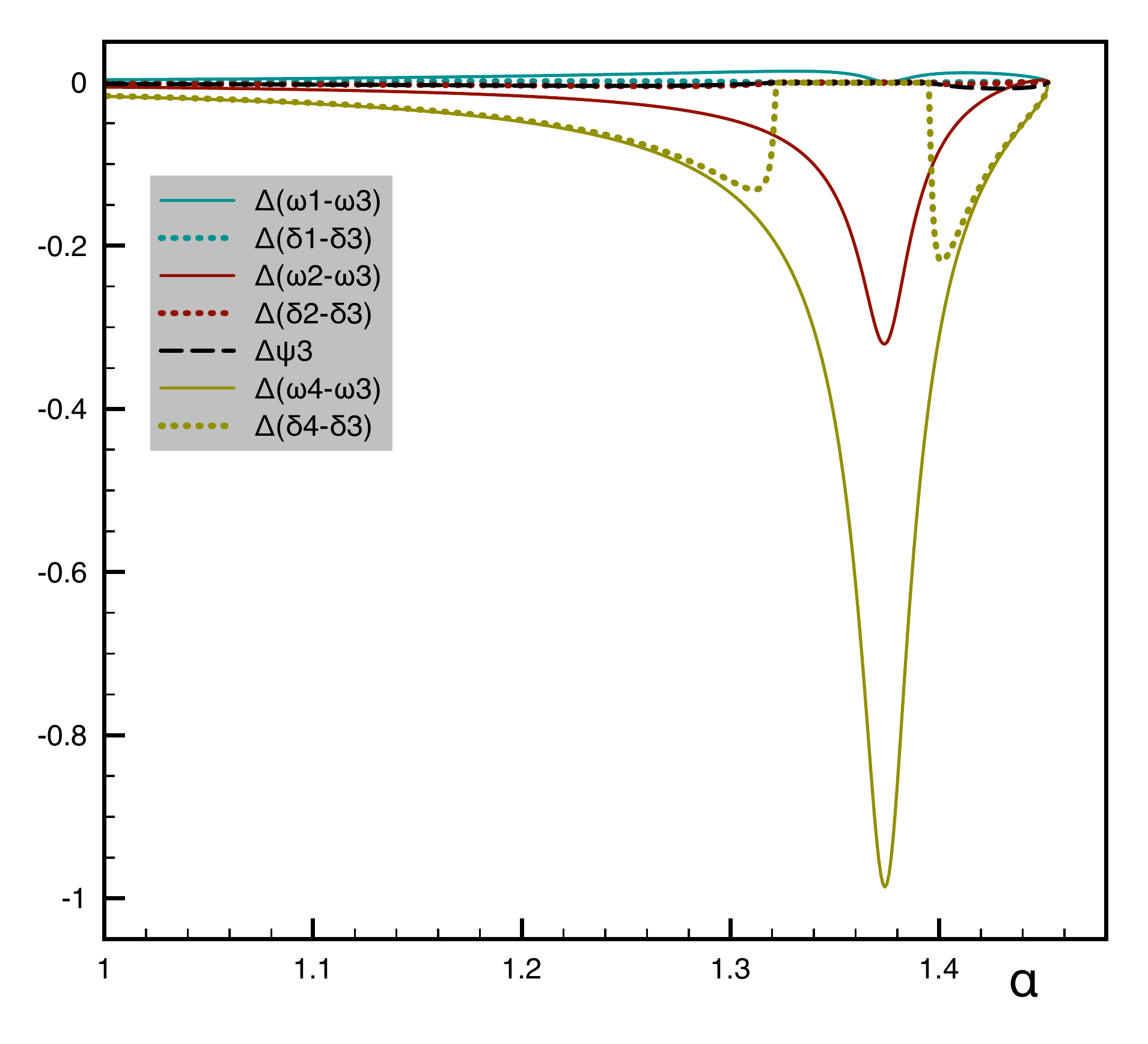}
   }
 \caption{Behavior of pair MISLPFs $\tilde{e}_{k(ij)}$ for the interaction of the oscillations $S_6$ and $S_7$ and different state variables.}
\label{fig:pair-1}
\end{figure} 

\begin{figure} [ht]
  \centerline{
  \includegraphics[width=7.8cm]{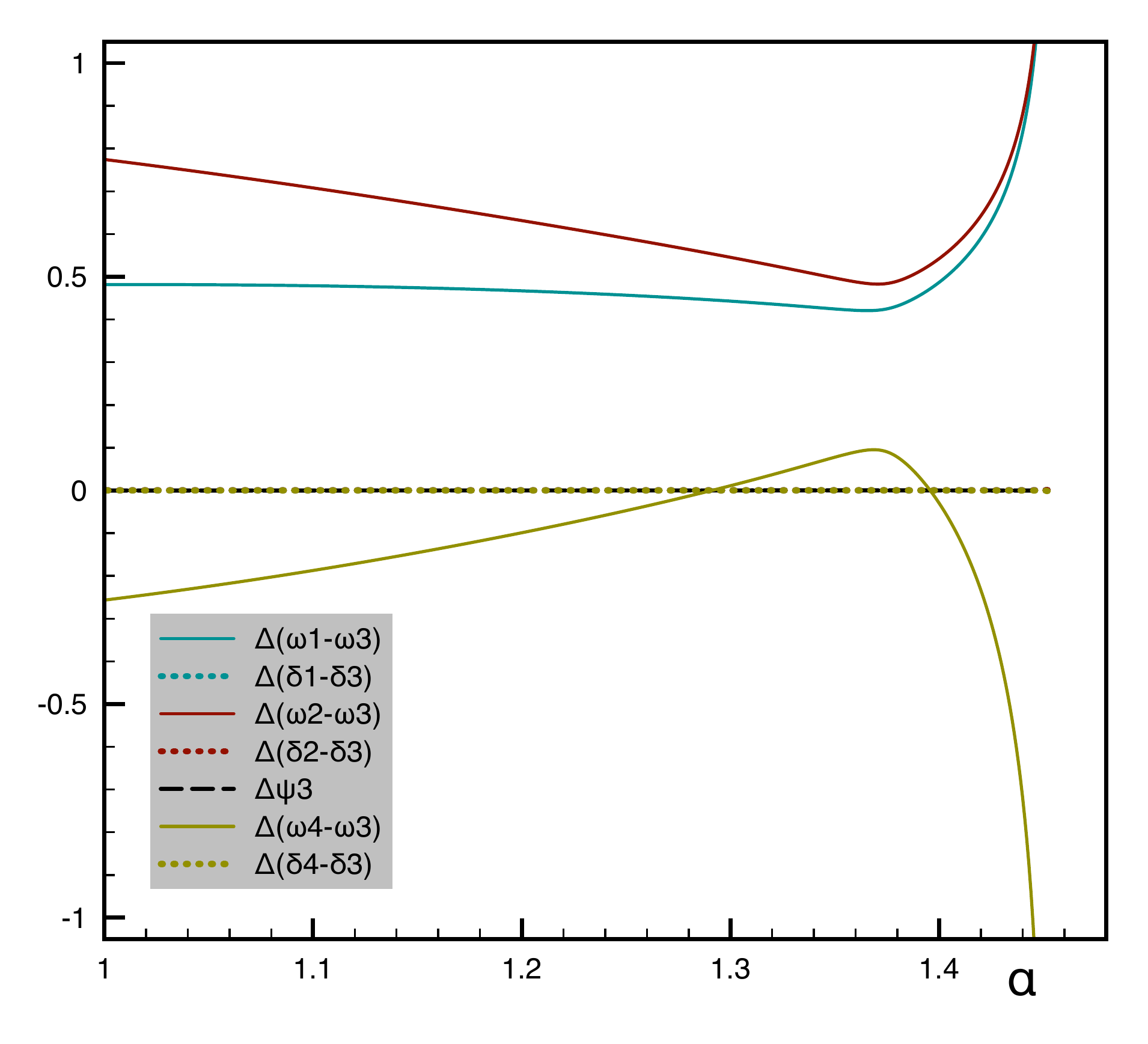}
   }
 \caption{Behavior of state participations in LMIEs $\overline{e}_{k(ij)}$ for the interaction of the oscillations $S_6$ and $S_7$ and different state variables.}
\label{fig:pair-2}
\end{figure} 

Figures~\ref{fig:pair-1} and \ref{fig:pair-2} show the behavior of pair MISLPFs $\tilde{e}_{k(ij)}$ defined in \eqref{def-5-1} and state participations in LMIEs $\overline{e}_{k(ij)}$ defined in \eqref{LMIE-pf}, respectively, for the interaction of the oscillations $S_6$ and $S_7$ and different state variables.
In Figure~\ref{fig:pair-1}, one can indicate the state variables $\Delta(\omega_2 - \omega_3)$ and $\Delta(\omega_4 - \omega_3)$, which are sensitive to this interaction, when it becomes resonant in the neighborhood of $\alpha_c \approx 1.375$. 
We note that the presence of the unstable mode $S_1$ in the system creates additional terms of the Lyapunov energy 
with large or infinite magnitude in some state variables (see Figure~\ref{fig:lyap-energies}).
This can make some state variables insensitive to the interaction of oscillations $S_6$ and $S_7$. 
In Figure~\ref{fig:pair-2}, the state variables $\Delta(\omega_1-\omega_3)$, $\Delta(\omega_2 - \omega_3)$, and $\Delta(\omega_4 - \omega_3)$, which provide the main contribution to this interaction, can be observed. Note that 
the magnitudes of $\overline{e}_{k(ij)}$ practically do not depend on the unstable mode $S_1$, 
as the Lyapunov energy of oscillations $S_6$ and $S_7$ does not depend on it.

\section{Conclusion}
This study proposes a novel LMA framework based on the concepts of Lyapunov energies and Lyapunov PFs, which characterize the time integrated energy associated either with particular modes and state variables, or with their pairwise combinations. It was proved that, in contrast to conventional PFs, the proposed indicators have characteristic properties that allow one to identify 
\begin{itemize}
\item the loss of stability of a particular mode, 
\item the resonant interactions between two modes, 
\item merging of two aperiodic modes into low-frequency oscillation, 
\end{itemize}
and associate these phenomena with certain system state variables. The calculation of the proposed Lyapunov indicators for the critical part of the spectrum does not require knowledge of the entire spectrum of the system matrix and can be performed independently. Therefore, the LMA can be performed quickly to analyze resonant interactions of the critical modes in large-scale dynamical systems. The proposed indicators can also be calculated using the sensitivities of eigenvalues obtained directly from measurements. 

Although the LMA was applied for analysis of the small-signal stability of the test power system in this work, its performance can also be tested for solving other problems of modal analysis, such as transient stability analysis, optimal placement of sensors and stabilizers, and cluster analysis of electrical networks, which is the subject of our further research.

\section*{Acknowledgment}
The authors thank Prof. P.Yu.~Chebotarev for providing helpful remarks. 
This work was supported by the Russian Science Foundation, project no.~19-19-00673.

\end{document}